\documentclass{amsart}
\usepackage{amssymb,bussproofs,pgf,float,listings,mdframed}
\usepackage{color,graphicx}
\usepackage[colorlinks]{hyperref}

\usepackage{tikz}
\usetikzlibrary{fit}
\usetikzlibrary{shapes}
\usetikzlibrary{arrows}
\usetikzlibrary{decorations.markings}
\usetikzlibrary{positioning}

\tikzset{mylabel/.style={font=\footnotesize}}
\tikzset{mylabelsmall/.style={font=\tiny}}
\tikzset{mymidlabel/.style={fill=white,font=\footnotesize}}
\tikzset{->-/.style={decoration={
      markings,
      mark=at position #1 with {\arrow{>}}},postaction={decorate}}}

\definecolor{mydark}{RGB}{73,68,62}
\definecolor{mymedium}{RGB}{128,129,135}
\definecolor{mylight}{RGB}{226,226,226}

\tikzstyle{smallcoqbox} = [fill=mylight!50, rectangle, inner sep=10pt, inner ysep=20pt]

\lstnewenvironment{coqcode}{\lstset{frame=single,backgroundcolor=\color{mylight},framerule=0pt}}{}

\newcommand{\UU}{\mathcal{U}}

\newcommand{\hprop}{\textnormal{\texttt{hProp}}}
\newcommand{\hset}{\textnormal{\texttt{hSet}}}

\newcommand{\iso}{\cong}

\newcommand{\id}[1]{\textnormal{\texttt{Id}}_{#1}}

\newcommand{\nat}{\mathtt{N}}
\newcommand{\zero}{\mathtt{0}}
\newcommand{\rr}{\mathtt{r}}
\newcommand{\successor}{\mathtt{S}}
\newcommand{\rec}{\textnormal{\texttt{rec}}}

\newcommand{\weq}{\texttt{WEq}}
\newcommand{\judge}[2]{(#1)\quad #2}

\newtheorem{theorem}{Theorem}[section]

\theoremstyle{definition}

\theoremstyle{remark}

\numberwithin{equation}{section}

\begin{document}
\newcounter{myi}
\setcounter{myi}{0}

\title[Homotopy type theory]{Homotopy type theory and
  Voevodsky's univalent foundations}


\author[\'{A}. Pelayo]{\'{A}lvaro Pelayo}
\address{Washington University\\  
  Mathematics Department \\
  One Brookings Drive, Campus Box 1146\\
  St Louis, MO 63130, USA, AND School of Mathematics,
  Institute for Advanced Study\\
  Einstein Drive, Princeton\\
  NJ 08540 USA}
\email{apelayo@math.wustl.edu, apelayo@math.ias.edu}
\thanks{Pelayo is partly supported by NSF CAREER Award DMS-1055897, Spain 
Ministry of Science Grant MTM 2010-21186-C02-01, and 
Spain Ministry of Science Sev-2011-0087.  Pelayo also received support
from NSF Grant DMS-0635607 during the preparation of this paper.}
\author[M. A. Warren]{Michael A. Warren}
\address{School of Mathematics,
  Institute for Advanced Study\\
  Einstein Drive, Princeton\\
  NJ 08540 USA.}
\curraddr{}
\email{mwarren@math.ias.edu}
\thanks{Warren is supported by the Oswald Veblen Fund and also
  received support from NSF Grant DMS-0635607 during the preparation
  of this paper.}


\date{\today}

\dedicatory{}

\begin{abstract}
 Recent discoveries have been made connecting
 abstract homotopy theory and the field of type theory from logic and
 theoretical computer science. This has given rise to a new field, which
 has been christened ``homotopy type theory''.  In this direction,
 Vladimir Voevodsky observed that it is possible to model type theory
 using simplicial sets and that this model satisfies an additional
 property, called the \emph{Univalence Axiom}, which has a number of
 striking consequences.  He has subsequently advocated a program, which
 he calls \emph{univalent foundations}, of developing mathematics
 in the setting of type theory with the Univalence Axiom and possibly
 other additional axioms motivated by the simplicial set model.  
 Because type theory possesses good computational properties, this
 program can be carried out in a computer proof assistant.  In this
 paper we give an introduction to homotopy type theory in Voevodsky's
 setting, paying attention to both theoretical and practical
 issues. In particular, the paper serves as an introduction to both
 the general ideas of homotopy type theory as well as to some of the
 concrete details of Voevodsky's work using the well-known proof
 assistant Coq.  The paper is written for a general audience of mathematicians
 with basic knowledge of algebraic topology; the paper does not assume any
 preliminary knowledge of type theory, logic, or computer science.
\end{abstract}

\maketitle


\section{Introduction}

Type theory is a branch of mathematical logic which developed out of
the work of Church \cite{Church:1933cl,Church:1940tu,Church:1941tc}
and which has subsequently found many
applications in theoretical computer science, especially in the theory
of programming languages \cite{Pierce:2002tp}.  For instance, the notion of
\emph{datatype} in programming languages derives from the type
theoretic notion of \emph{type}.  Recently, a number of deep and
unexpected connections between a form of type theory (introduced
by Martin-L\"{o}f
\cite{MartinLof:1998tw,MartinLof:1975tb,MartinLof:1982bn,MartinLof:1984tr}) and
homotopy theory have been discovered, opening the way to a new area of research in
mathematics and theoretical computer science which has recently been
christened \emph{homotopy type theory}.  Due to the nature of the mathematical
results in this area, we believe that there is great potential for the
future research in this area to have a considerable impact on a number
of areas of pure and applied mathematics, as well as on the practice
of mathematicians.

In 1998, Hofmann and Streicher \cite{Hofmann:1998ty} constructed a model
of Martin-L\"{o}f type theory in the category of
groupoids.  They also observed that the data of type theory
itself naturally gives rise to a kind of $\infty$-groupoid structure (although they did not prove this for any precise
definition of $\infty$-groupoid).  In 2001, Moerdijk
speculated that there should be some connection between Quillen model
categories and type theory.  Then between 2005 and 2006 Awodey
and Warren \cite{Awodey:2009bz,Warren:2006vf,Warren:2008ts}, and
Voevodsky \cite{Voevodsky:VSNHLC,Voevodsky:2009,Vo2012a,Vo2012} independently understood how to
interpret type theory using ideas from homotopy theory (in the
former case, using the general machinery of Quillen model categories
and weak factorization systems, and in the latter case using
simplicial sets).  Subsequently, around 2009, Voevodsky
\cite{Voevodsky:2009} realized that
the model of type theory in simplicial sets satisfies an additional
axiom, which he called the \emph{Univalence Axiom}, that is not in general
satisfied.  Crucially, satisfaction of the Univalence Axiom is a
property which distinguishes Voevodsky's model of type theory in
simplicial sets from the more familiar set theoretic model (it does
not hold in the latter).

These results and others (described in more detail below) give rise to
what might be called the \emph{univalent perspective}, wherein one
works with a kind of type theory augmented by additional axioms, such as the
Univalence Axiom, which are valid in the simplicial set model.  In
this approach one thinks of types as spaces or homotopy
types and, crucially, one is able to manipulate spaces directly
without having first to develop point set topology.\footnote{One
  should be careful here since the spaces
  under consideration should be suitably nice spaces (fibrant and
  cofibrant spaces in the language of homotopy theory) on the one
  hand.  On the other hand, when thinking of homotopy types one should
  not think of homotopy types living in the homotopy category, but
  rather homotopy types living (via their presentations as fibrant and
  cofibrant spaces) in the category of spaces (simplicial sets).}  Although
interesting in its own right, this perspective becomes significantly
more notable in light of the good computational properties of the
kind of type theory employed here.  In particular, type theory of
the sort considered here forms the underlying theoretical framework
of several computer proof assistants such as \emph{Agda} and
\emph{Coq} (see \cite{Coquand} and \cite{Bertot:2004uj},
respectively).  

In practical terms, this means that it is possible to
develop mathematics involving spaces in computer systems which are
capable of verifying the correctness of proofs and of providing some
degree of automation of proofs.  We refer the reader to
\cite{Simpson:2004bt} and \cite{Hales:2008ud} for two accounts of
computer proof assistants (and related developments) written for a
general mathematical audience.  

Voevodsky has, since his discoveries mentioned above, been advocating
the formalization of mathematics in proof assistants, as well as
greater interaction between the developers of computer proof
assistants and pure mathematicians. He has himself written thousands
of lines of code in the Coq proof assistant, documenting topics
ranging from the development of homotopy theoretic notions and proofs
of new results in type theory, to the formalization of the basics of
abstract algebra.

Voevodsky's univalent perspective, as detailed in his Coq files, is a
unique view of mathematics and we believe that it deserves to be more widely known.
Unfortunately, for a mathematician without some background in type
theory, homotopy theory and category theory, we believe that the
prospect of reading thousands of lines of Coq code is likely rather
daunting (indeed, it may be daunting even for those with the
prerequisites listed above).  In this paper we attempt to
remedy this by providing an introduction to both homotopy type theory
and the univalent perspective, as well as to some of the material contained
in Voevodsky's Coq files.  It is our hope that the reader who is not
interested in the Coq code, but who is curious about homotopy type
theory will benefit from an account of this field specifically
targeted at a general mathematical audience.  For those who are
interested in the Coq code, we believe that this paper can act as an
accessible introduction.  Indeed, it is our ultimate aim that this
paper will encourage other mathematicians to become involved
in this area and in the use of computer proof assistants in general.
This paper also serves as a guide to the  authors' recent on-going
work on $p$-adic arithmetic and $p$-adic integrable systems in Coq
\cite{PeVoWa2012}.

We believe that the timing of this article is perspicuous in part
because there is a Special Year devoted to Voevodsky's program
during the 2012-2013 academic year at the Institute for Advanced
Study. Further information on the activities of this program is
given in  Awodey, Pelayo and Warren \cite{AwPeWa2012}.

\subsubsection*{Disclaimers}

This article is aimed at mathematicians who want to understand the basics
of homotopy type theory and the univalent perspective. It is
written for a broad audience of readers who are not necessarily
familiar with type theory and homotopy theory. Because of the
introductory nature of the article we are less precise than one would
be in a research article.  This is especially true when it comes to
describing type theory and Coq, where we eschew excessive terminology
and notation in favor of a more informal approach.  For those readers
with the requisite background in logic and category theory who are
interested in a more detailed account we refer to \cite{Awodey:TTH}.  Needless
to say, the present article has no intention of being comprehensive,
it is merely an invitation to a new and exciting subject.

It is worth mentioning that there are already a number of introductions
to Coq available, see for instance \cite{Bertot:2004uj}, which are far
more comprehensive and precise than this article in their treatment of
the proof assistant itself.  However, such introductions inevitably
make use of features of the Coq system which do not enter into (or are
even incompatible with) the univalent perspective.  As such, we warn
the reader that this paper is \emph{not} a Coq tutorial: it is an
introduction to the univalent perspective which along the way
also describes some of the basic features of Coq.

Finally, Coq is not an automatic theorem prover, but rather an interactive
theorem prover: it helps one to verify the correctness of proofs which
are themselves provided by the user.\footnote{Technically, it \emph{is}
  possible to automate proofs to a large extent in Coq (via the built-in
\emph{tactics language}), but, aside
  from a minor amount of automation implicit in the ``tactics'' we
  employ below, we will not go into details regarding these features
  of Coq.  The interested reader might consult \cite{Chlipala:CP} for a good
  introduction to Coq which pays particular attention to automation.}

\section{Origins, basic aspects and current research} \label{sec:origins}

This section gives an overview of the origins of homotopy type
theory and Voevodsky's univalent perspective.  We assume that
readers are more familiar with algebraic topology than with type
theory (indeed, we expect that many readers may have not have been
exposed to \emph{any} type theory).  As such, we will introduce type
theory ---already with the homotopy theoretic interpretation of \cite{Awodey:2009bz,Voevodsky:2009} in mind--- by analogy with certain developments and constructions
in algebraic topology.  

This approach is admittedly anachronistic, but
we hope that it will serve as an accessible starting point for readers
coming from outside of type theory.  Along the way we will try to give
an idea of the historical development of the field.  However, we make
no attempt to provide a comprehensive history of homotopy theory or type theory.
For the early history of homotopy theory we refer to
\cite{Dieudonne:2009dv}.

Before starting with type theory it is perhaps worth remarking that
type theory (like set theory) is a logical theory which is given by a
collection of \emph{rules}.  Anyone interested in type theory
should at some point study these rules, but doing so is not
strictly necessary in order to give some flavor of the theory.  As
such, we choose to abstain from giving a formal presentation of type
theory.

\subsection{The homotopy theoretic interpretation of type theory}\label{sec:homotopy_interp}

Although the mathematical notion of \emph{type} first appears in
Russell's \cite{Russell:1903wn} work on the
foundations of mathematics, it was not until the work of Church
\cite{Church:1940tu} that type theory in its modern form was born.
Later, building on work by Curry \cite{Curry:1934vy}, Howard \cite{Howard:FTNC}, Tait
\cite{Tait}, Lawvere \cite{Lawvere:2006tl},
Scott \cite{Scott:1970vu} and others, Martin-L\"{o}f
\cite{MartinLof:1998tw,MartinLof:1975tb,MartinLof:1982bn,MartinLof:1984tr}
developed a generalization of Church's system which is now called
\emph{dependent} or \emph{Martin-L\"{o}f type theory}. We will be
exclusively concerned with this form of type theory and so the
locution ``type theory'' henceforth refers to this particular
system.\footnote{We caution the reader that there are indeed many
  different kinds of type theory.  Explicitly, we are concerned with
  the \emph{intensional} form of dependent type theory.}

The principal form of expression in type theory is the statement that
the \emph{term} $a$ is of
\emph{type} $A$, which is written as 
\begin{align*}
  a:A.
\end{align*}
There are a number of ways that the expression $a:A$ has traditionally
been motivated:
\begin{enumerate}
\item $A$ is a set and $a$ is an element of $A$.
\item $A$ is a problem and $a$ is a solution of $A$.
\item $A$ is a proposition and $a$ is a proof of $A$.
\end{enumerate}
Roughly, of the perspectives enumerated here,
the first is due to Russell \cite{Russell:1903wn}, the
second is due to Kolmogorov \cite{Kolmogoroff:1932wl}, and the third
---usually called the \emph{Curry-Howard correspondence}---
is due to Curry and Howard \cite{Howard:FTNC}.  We will say more about
these three motivations in the sequel.

The starting point for understanding the connections between homotopy
theory and type theory is to consider a fourth alternative to these
motivations:
\begin{enumerate}
\item[(4)] $A$ is a space and $a$ is a point of $A$.
\end{enumerate}
We will be intentionally vague about exactly what kinds of
spaces we are considering, but we recommend that readers have in mind,
e.g., topological spaces, (better yet) CW-complexes, or (still better yet)
Kan complexes (in which case, ``point'' means $0$-simplex).  The
remarkable thing about (4) is that it helps to clarify certain features of
type theory which originally seemed odd, or
even undesirable, from the point of view of the
interpretations (1) -- (3) above.  

In addition to the kinds of types and terms described above, we also
may consider types and terms with parameters.  These are usually called
\emph{dependent} types and terms.  E.g., when $B$ is a type, we might
have a type
\begin{align*}
  \judge{x:B}{E(x),}
\end{align*}
which is parameterized by $B$ (here $x$ is a variable).  In terms of the motivation (1) above
in terms of sets, we would think of this as a $B$-indexed family $(E_{x})_{x\in
  B}$.  From the homotopy theoretic point of view (4), we think of
such a type as describing a fibration $E\to B$ over the space $B$.
Similarly, we think of a parameterized term
\begin{align*}
  \judge{x:B}{s(x):E(x)}
\end{align*}
as a continuous section of the fibration $E\to B$.
\begin{figure}[H]
  \begin{tikzpicture}[color=mydark, fill=mylight, line width=1pt]
    \filldraw (0,0) ellipse (1 and .5);
    \node (bb) at (1.75,0) {$\;B.$};
    \filldraw (0,1.5) ellipse (1 and .5);
    \filldraw[color=mylight,opacity=5] (-1,1.5) rectangle (1,2.75);
    \draw[inner sep=0] (-1,1.48) to (-1,2.75);
    \draw (1,1.48) to (1,2.75);
    \draw[opacity=.5] (1,1.5) arc (0:180:1 and .5);
    \node[circle,fill=mydark,inner sep=0pt,minimum size=.75mm] at
    (-.5,0) {};
    \node[font=\tiny] at (-.35,.05) {$x$};
    \draw[opacity=.6] (-.5,1.45) to (-.5,2.75);
    \node[circle,fill=mydark,inner sep=0pt,minimum size=.75mm] at
    (-.5,1.7) {};
    \node[font=\tiny] at (-.2,1.7) {$s(x)$};
    
    \filldraw (0,2.75) ellipse (1 and .5);
    \draw[opacity=.5] (-.5,1.45) to (-.5,2.8);
    
    \node (ee) at (1.75,2.25) {$E$};
    \draw[->] (ee) to node[auto,swap] {$p$} (bb);
    \draw[->,bend right] (bb) to node[auto,swap] {$s$} (ee);
  \end{tikzpicture}

  \begin{tabular}[h!]{lr}
    \emph{type theory} &  \emph{homotopy theory}\\
    \hline
    $\judge{x:B}{E(x)}$ & $p:E\to B$ is a fibration over $B$\\
    $\judge{x:B}{s(x):E(x)}$ & $s$ is a section of $p$
  \end{tabular}
  \caption{Homotopy theoretic interpretation of dependent types and terms.}
\end{figure}
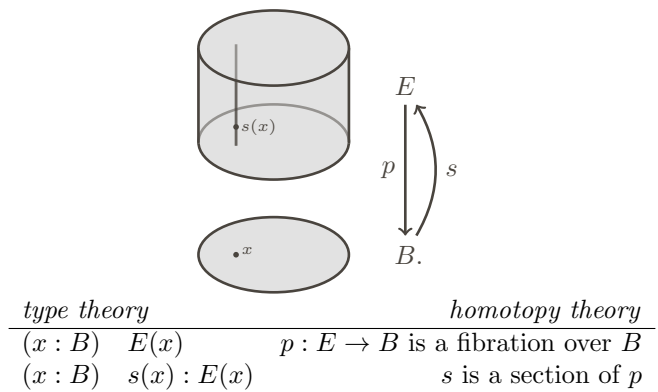

Of course, none of this would be useful without being given some
types and terms to start with and some rules for generating new types and terms
from old ones.  It is to these that we now turn.

\subsection{Inductive types}\label{sec:inductive_intro}

Among all of the types which can be constructed, some of the most
significant and interesting are the \emph{inductive types}.  The most
familiar example of an inductive type is the type $\nat$ of natural
numbers.  For types, the property of being inductive is the same as
being free in an appropriate (type theoretic) sense.  E.g., $\nat$ is
freely generated by the following data:
\begin{itemize}
\item a term $\zero:\nat$; and
\item a term $\judge{x:\nat}{\successor(x):\nat}$.
\end{itemize}
Consequently, $\nat$ has a type theoretic universal property which
should be familiar as proof by induction (alternatively, definition by
recursion).  Namely, given any type 
\begin{align*}
  \judge{x:\nat}{E(x)}
\end{align*}
fibered over $\nat$ together with terms:
\begin{itemize}
\item $e:E(\zero)$ (\emph{base case}); and
\item $\judge{x:\nat,y:E(x)}{f(x,y):E(\successor(x))}$
  (\emph{induction step}),
\end{itemize}
there exists a term, which we denote by
$\judge{x:\nat}{\rec(e,f,x)}$, of type $E(x)$.  This term has
the corresponding properties that 
\begin{align*}
  \rec(e,f,\zero) & = e,\text{ and}\\
  \rec(e,f,\successor(x)) & = f(x,\rec(e,f,x)).
\end{align*}
This should be compared with the usual way that functions from the
natural numbers are defined by recursion.

In general, it is possible to construct inductive types which are
generated (in the sense that $\nat$ is freely generated by zero and
successor) by arbitrary generators (subject to some technical
conditions on the kinds of generators allowed).  Indeed, it is also
possible to construct inductive types which
are themselves dependent.  We will encounter several other examples of inductive
types below, but for now we mention one single significant example.

If $A$ is a type, consider the inductive type fibered over $A\times A$
with a single generator $\rr(a)$ in the fiber over the pair $(a,a)$.
That is, consider the smallest (in an appropriate sense) fibration
over $A\times A$ having such elements $\rr(a)$ in the fibers.
Somewhat miraculously, this inductive type turns out to be a well
known topological space: the space of all (continuous) paths in $A$.
\begin{theorem}[Awodey and Warren \cite{Awodey:2009bz}]\label{theorem:aw}
  The path space fibration $A^{I}\to A\times A$ is an inductive
  type.
\end{theorem}
Actually, Theorem \ref{theorem:aw} is just a small part of a more
general result from \emph{ibid} (see also \cite{Warren:2008ts} for
further details): it is possible to model type theory in
weak factorization systems or Quillen model categories
\cite{Quillen:1967uz} which satisfy certain coherence conditions.
Going the other way, Gambino and Garner \cite{Gambino:2008bg} showed
that it is possible to construct a weak factorization system from the
syntax of type theory.

\subsection{Groupoids and $\infty$-groupoids}\label{sec:groupoids}

The type theoretic content of Theorem \ref{theorem:aw} is more than it
perhaps looks at first.  The reason for this is that the inductive
type corresponding to the path space fibration is in fact a well-known
inductive type and one which plays a crucial role in type theory.  Called the
\emph{identity type of $A$}, this type is usually written
type-theoretically as 
\begin{align*}
  \judge{x:A,y:A}{\id{A}(x,y)}
\end{align*}
with generators written as
\begin{align*}
  \judge{x:A}{\rr(x):\id{A}(x,x)}.
\end{align*}
Type theorists traditionally thought of the type $\id{A}(a,b)$ as
something like the ``proposition that $a$ and $b$ are equal proofs
that $A$''.  It was also known that it is possible to
construct a set theoretic model of type theory by thinking, as in
perspective (1) above, of types as sets and terms as elements.  In
this set theoretic interpretation we have that
\begin{align*}
  \id{A}(a,b) & =
  \begin{cases}
    1 & \text{ if }a=b,\text{ and}\\
    0 & \text{ otherwise.}
  \end{cases}
\end{align*}
I.e., $\id{A}(a,b)$ is either the one point set or the empty set
depending on whether or not $a$ and $b$ are in fact equal.  Under this
interpretation it follows that if there is a term $f:\id{A}(a,b)$,
then in fact $a=b$.  This was originally seen by many, no doubt based
on intuitions gleaned from the set theoretic model, as a desirable
property and this property was added by many type theorists (including, for a time,
Martin-L\"{o}f himself) as an axiom, which we call the
\emph{0-truncation rule}.  

Type theory with the 0-truncation rule may
be somewhat easier to work with and, by work of Seely
\cite{Seely:1984jw} and a coherence result due to Hofmann \cite{Hofmann:1995be}, it is
possible, using the machinery of locally cartesian closed categories,
to obtain many different models of type theory which satisfy the
0-truncation rule.  Nonetheless, adding the 0-truncation rule destroys
the good computational properties of type theory (see Section
\ref{section:computational} for a brief description of the
computational aspects of type theory).  Many of the early questions in
Martin-L\"{o}f type theory were related to questions about the
behavior of the identity types in the theory without 0-truncation.  In
particular, many results were concerned with trying to understand
whether facts which are consequences of the 0-truncation axiom also
hold without the 0-truncation axiom.  Streicher's
\emph{Habilitationsschrift} \cite{Streicher:vg} contains many fundamental results in this
direction.

One important example of the kind of questions about identity types
that type theorists were trying to answer is the problem of
``uniqueness of identity proofs'': if $f$ and $g$ are both of type
$\id{A}(a,b)$, does it follow that $f=g$ (or even that there exists a
term $\eta:\id{\id{A}(a,b)}(f,g)$)?  In order to solve this problem
Hofmann and Streicher \cite{Hofmann:1998ty} constructed a model of
type theory in which types are interpreted as groupoids and fibrations
of groupoids, and terms are interpreted as sections (see below for the
notion of groupoids).  In the process of
constructing this model, Hofmann and Streicher discovered an
interesting fact which we will now describe.

Given terms $a$ and $b$ of type $A$, there is an equivalence relation
$\simeq$ on the set of terms of type $\id{A}(a,b)$ given by setting
$f\simeq g$ if and only if there exists a term
$\eta:\id{\id{A}(a,b)}(f,g)$.  In terms of spaces, $f$ and $g$
correspond to paths from $a$ to $b$ in the space $A$ and $\eta$
corresponds to a homotopy rel endpoints from $f$ to $g$.  Hofmann and
Streicher realized that the quotient of the set
of terms $f:\id{A}(a,a)$ modulo $\simeq$ forms a group.  In fact, they
realized that the type $A$ can be made into a groupoid.  Recall that a
groupoid is a category in which every arrow is invertible.  To turn
$A$ into a groupoid we take the objects to be the terms $a:A$ and the
edges to be equivalence classes of $f:\id{A}(a,b)$ modulo $\simeq$.
We now see that these two constructions correspond, under the homotopy
theoretic interpretation of type theory sketched above, to the
constructions of the fundamental group $\pi_{1}(A,a)$ of the space $A$
with basepoint $a$ and the fundamental groupoid $\Pi_{1}(A)$ (see
\cite{Brown:2006tj} for more on the fundamental groupoid of a space).
In fact, Hofmann and Streicher realized that their construction seemed to give
some kind of $\infty$-groupoid, but they did not pursue
this possibility.  The first non-syntactic higher-dimensional models (which
are shown to satisfy all of the required coherence conditions) of type theory
appeared later in work of Voevodsky \cite{Voevodsky:2009,Kapulkin:USS} and Warren
\cite{Warren:2008ts,Warren:2011tn}.

In homotopy theory one is also concerned with groupoids and
$\infty$-groupoids.  The fundamental groupoid $\Pi_{1}(A)$
of a space $A$ is the basepoint-free generalization of Poincar\'{e}'s
\cite{Poincare:1895vx} fundamental group and captures the homotopy
1-types.  Here a \emph{homotopy $n$-type}
is intuitively a space $A$ for which the higher-homotopy groups
$\pi_{i}(A,a)$, for $i>n$, vanish.  For spaces $A$ that are not
1-types one must consider higher-dimensional generalizations of
the notion of groupoid in order to capture the homotopy theoretic
content of $A$.  In his letter to Quillen, Grothendieck \cite{Grothendieck:2005wy} emphasized the importance of
finding an infinite dimensional generalization of the notion of
groupoid which would capture the homotopy types of spaces
and indeed he offered a suggestion himself (see \cite{Maltsiniotis}
for a modern exposition of Grothendieck's definition).  This problem has been one of
the leading motivations for the development of higher-dimensional
category theory (see \cite{Baues:1995js} for a
detailed overview of the problem of modeling homotopy types, and see
\cite{Loday:1982vg,Kapranov:1991uw,Simpson:1998wi,Tamsamani:1999ia,Cisinski:2007by}
for some results on modeling homotopy types).

In the setting of type theory, with respect to the Batanin-Leinster
\cite{Batanin:1998kp,Leinster:2004fi} notion of $\infty$-groupoid,
types have associated $\infty$-groupoids:
\begin{theorem}[van den Berg and Garner \cite{vandenBerg:2011ec}, and
  Lumsdaine \cite{Lumsdaine:2010ew}]
  Every type $A$ has an associated fundamental $\infty$-groupoid $\Pi_{\infty}(A)$.
\end{theorem}
Hence the set theoretic intuition for
the meaning of types fails to accurately capture certain features of
the syntax, whereas those features (non-trivial higher-dimensional
structure) are captured by the homotopy theoretic interpretation of
types.

\subsection{The univalent model of type theory}\label{sec:univalent_model}

We mentioned above the problem of finding a notion of
$\infty$-groupoid which completely captures the notion of homotopy
type.  One such notion is provided by \emph{Kan complexes}.
Introduced by Kan \cite{Kan:1957up}, these were shown by Quillen
\cite{Quillen:1967uz} to provide a model of homotopy types.  Kan
complexes are simplicial sets which satisfy a certain combinatorial
condition.  In the work of Joyal \cite{Joyal:2008uc} and Lurie
\cite{Lurie:2009un} on $\infty$-toposes, the Kan complexes are the
$\infty$-groupoids.  The starting place for Voevodsky's univalent
perspective is the following result:
\begin{theorem}[Voevodsky \cite{Voevodsky:2009}]
  Assuming the existence of Grothendieck universes (sufficiently large
  cardinals), there is a model of type theory in the category of
  simplicial sets (equipped with well-orderings) in which types are interpreted as Kan complexes and
  Kan fibrations.
\end{theorem}
We will henceforth refer to Voevodsky's model as the \emph{univalent
  model} of type theory.  The particular kind of type theory
considered by Voevodsky includes a
type $\UU$ which is a universe of types which we refer to  as
\emph{small types}.\footnote{The universes
  in Voevodsky's model play two roles.  The first is simply to
  interpret universes of types.  The second is to avoid certain
  coherence issues which arise in the homotopy theoretic
  interpretation of type theory.}  Given small types $A$ and $B$, it
is then natural to ask what kinds of terms arise in the identity type
$\id{\UU}(A,B)$.  Voevodsky realized
that, although this type \emph{a priori} possesses no interesting
structure, in the univalent model it turns out to be non-trivial.
Based on this observation, Voevodsky proposed to add to the axioms of
type theory an additional axiom, called the \emph{Univalence Axiom},
which would ensure that the identity type of $\UU$ behaves as it does
in the univalent model.  We will now explain this axiom and some of
its consequences.

It will be instructive to compare several different ways of
understanding the Univalence Axiom.  However, we will first start by
giving an explicit description of the axiom.  Given small types $A$
and $B$ there is, in addition to the identity
type $\id{\UU}(A,B)$, a type $\weq(A,B)$ of \emph{weak equivalences}
from $A$ to $B$.  Intuitively, thinking of $A$ and $B$ as spaces, a
weak equivalence $f:A\to B$ is a continuous function which induces
isomorphisms on homotopy groups: 
\begin{align*}
  \pi_{n}(f):\pi_{n}(A,a)\iso\pi_{n}(B,f(a))
\end{align*}
for $n\geq 0$.\footnote{Technically, because we are dealing with
  sufficiently nice spaces such as CW-complexes or Kan complexes, the
  weak equivalences do in fact coincide with the homotopy equivalences.}
Since the identity $1_{A}:A\to A$ is a weak equivalence, there is, by
the induction principle for identity types, an induced map
$\iota:\id{\UU}(A,B)\to\weq(A,B)$ and the Univalence Axiom can be
stated as follows:
\begin{description}
\item[Univalence Axiom \textnormal{(Voevodsky)}] The map $\iota:\id{\UU}(A,B)\to\weq(A,B)$
  is a weak equivalence.
\end{description}
That is, the Univalence Axiom imposes the condition that the identity
type between two types is naturally weakly equivalent to the type of
weak equivalences between these types.  The Univalence Axiom makes it
possible to automatically transport constructions and proofs between
types which are connected by appropriately defined weak
equivalences.

Without the Univalence Axiom there are three \emph{a priori} different ways in which two small types
$A$ and $B$ can be said to be equivalent:
\begin{enumerate}
\item $A=B$.
\item There exists a term $f:\id{\UU}(A,B)$.
\item There exists a term $f:\weq(A,B)$.
\end{enumerate}
The Univalence Axiom should be understood as asserting (in a type
theoretic way) that (2) and (3) coincide.  (Interestingly, the
0-truncation axiom for $\UU$ asserts that (1) and (2)
coincide.)  That is, the Univalence Axiom answers the question ``What is a path from $A$ to $B$ in
the space of small spaces?'' by stipulating that such a path
corresponds to a weak equivalence $A\to B$.  

Alternatively, the Univalence Axiom may be understood as stating that the types
of the form $\weq(A,B)$ are inductively generated by the identity maps
$1_{A}:\weq(A,A)$.  Part of the appeal of the Univalence Axiom is that
it has a number of interesting consequences which we will discuss
below.  The connection between the Univalence Axiom and object
classifiers from topos theory has recently been investigated by
Moerdijk \cite{Moerdijk:2012vl}.

\subsection{The univalent perspective}

Following Voevodsky, we define a filtration of types by what are
called \emph{h-levels} extending the usual
filtration of spaces by homotopy $n$-types.  The h-levels are defined
as follows:
\begin{itemize}
\item A type $A$ is of h-level $0$ if it is contractible.
\item A type $A$ is of h-level $(n+1)$ if, for all terms $a$ and $b$
  of type $A$, the type $\id{A}(a,b)$ is of h-level $n$.
\end{itemize}
For $n>1$, $A$ is of h-level $n$ if and only if it is a homotopy
$(n-2)$-type.  E.g., types of h-level $2$ are the same as homotopy
0-types: spaces which are homotopy equivalent to sets.  From a more
category theoretic point of view we can view the h-levels as
follows:
\begin{figure}[H]
  \centering
  \begin{tabular}{cl}
    \emph{h-level} & \emph{corresponding spaces up to weak equivalence}\\
    \hline
    0 & the contractible space $1$\\
    1 & the space $1$ and the empty space $0$\\
    2 & the homotopy $1$-types (i.e., groupoids)\\
    3 & the homotopy $2$-types (i.e., 2-groupoids)\\
    \vdots & \vdots\\
    $n$ & the homotopy $(n-2)$-types (i.e., weak $(n-2)$-groupoids)\\
    \vdots & \vdots
  \end{tabular}
  \caption{h-levels.}
  \label{fig:hlevels}
\end{figure}

It is worth recording several basic observations.  First, weak
equivalences respect h-level: if there is
a weak equivalence $f:A\to B$, then $A$ is of h-level $n$ if and only
if $B$ is of h-level $n$.  Secondly, h-levels are cumulative in the sense
that if $A$ is of h-level $n$, then it is also of h-level $(n+1)$.
Finally, for any $n$, if $A$ (or $B$) is of h-level $(n+1)$ then so is
$\weq(A,B)$.

We denote by $\hprop$ the type of all (small) types of h-level $1$.
The type $\hprop$ plays the same role, from the univalent perspective,
as the boolean algebra $\mathbf{2}:=\{0,1\}$ in classical logic and set
theory, or the subobject classifier $\Omega$ in topos theory.  We will
usually refer to types in $\hprop$ as \emph{propositions}.  For
$\hprop$, the Univalence Axiom states that paths in $\hprop$
correspond to logical equivalences.  I.e., for propositions $P$ and
$Q$, it is a necessary and sufficient condition for there to exist a
term of type $\id{\hprop}(P,Q)$ that $P$ and $Q$ should be logically
equivalent.

Similarly, we denote by $\hset$ the type of all (small) types of
h-level $2$ and we refer to these types as \emph{sets}.  The type
$\hprop$ is itself a set.  To see this, note that, by the
``propositional'' form of the Univalence Axiom mentioned above, there is
a weak equivalence $\iota:\id{\hprop}(A,B)\to\weq(A,B)$.  It then
follows from the basic observations on h-levels summarized above that
$\id{\hprop}(A,B)$ has h-level 1, as required.  This result is a
special case of the more general fact that the Univalence
Axiom implies that the type $\textnormal{\texttt{hlevel}}_{n}$ of all
(small) types of h-level $n$ is itself of h-level $(n+1)$.

As we see it, the principal idea underlying the univalent perspective is that,
rather than developing mathematics in the setting of set theory where one must
``build'' all of mathematics up from the emptyset and the operations
of set theory, we should instead work in a formal system (namely, type
theory) where we are given at the outset the world of spaces (homotopy types).
In this setting we would still have all of the sets available to us,
but they are ``carved out of'' or extracted from the
universe as the types of h-level 2.

\begin{figure}[H]
  \centering
  \begin{tabular}{l|cc}
    & \emph{set theoretic} & \emph{univalent}\\
    \hline
    \emph{spaces} & constructed & given\\
    \emph{sets} & given & extracted
  \end{tabular}
  \caption{Sets and spaces from set theoretic and univalent perspectives.}
  \label{fig:comparison}
\end{figure}
Something which is not revealed in this simple comparison is that it
is considerably easier to extract sets from the world of homotopy
types than it is to construct spaces from sets.

From the univalent perspective, the development of ordinary ``set-level
mathematics'', which deals with sets and structures (e.g., groups,
rings, \ldots) on sets, is quite similar with the ordinary development
of mathematics.  However, in this setting it is also easy to develop
``higher-level mathematics''.  To take on simple example, the notion of
\emph{monoid} can be axiomatized in the ordinary way.  I.e., a monoid
consists of a type $M$ together with a binary operation $\mu:M\times
M\to M$ which is associative and unital (in the appropriate type
theoretic sense).  When we restrict $M$ to just
small types in $\hset$, we obtain the usual notion of monoid.
However, when $M$ is allowed to be an arbitrary type, we obtain the
notion of \emph{homotopy associative H-space} in the sense of Serre
\cite{Serre:1951uf}.  We believe that it is an advantage of the
univalent perspective that it is in fact easier to work with such
higher-level structures in the univalent setting than in the familiar
set theoretic setting.

In addition to the fact that it is efficient to reason about
higher-dimensional structures in the univalent setting, there are also
technical advantages to doing so.  For example, in the presence of the
Univalence Axiom, any structure on a type $A$ which is type
theoretically definable can be transferred along a weak equivalence
$A\to B$ to give a corresponding structure on $B$.  In general, being
able to transfer structures along weak equivalences (or even homotopy
equivalences) is non-trivial (see, e.g., \cite{Loday:2012vi} for some
examples of such ``homotopy transfer theorems'' and their
consequences).  Therefore, being able to work in a setting where such
transfer is ``automatic'' is technically quite convenient.

\subsection{Computational aspects}\label{section:computational}

One of the advantages of working with the particular flavor of
Martin-L\"{o}f type theory employed in the univalent setting is that
this theory has good computational properties.  In technical terms,
type checking in this theory is decidable.  Consequently, it is
possible to implement the theory on a computer.  This is essentially
what has been done in the case of the ``proof assistants'' Coq and
Agda. Therefore, mathematics in the univalent setting can be
formalized in these systems and the veracity of proofs can be
automatically checked.  In the case of reasoning involving homotopy
theoretic or higher-dimensional algebraic structures, which sometimes
involve keeping track of large quantities of complex combinatorial
data (think of, e.g., reasoning involving \emph{tricategories}
\cite{Gordon:CT}), being able to make use of the computer to ensure
that calculational errors have not been made is potentially quite
useful.

Part of the reason that Martin-L\"{o}f type theory enjoys such good
computational properties is that it is a \emph{constructive} theory.
\emph{Classical logic} is the usual logic (or framework for
organizing mathematical arguments) employed in mathematics (it is the
logic of the Boolean algebra $\{0,1\}$).  The logic employed in
constructive mathematics is obtained from classical logic by omitting
the \emph{law of excluded middle}, which stipulates that, for any
statement $\varphi$, either $\varphi$ or not $\varphi$.  Working
constructively is often more challenging than working classically and
sometimes leads to new developments.  Although there are many reasons that one
might be interested in pursuing constructive mathematics we will give
several practical reasons.  First, constructive mathematics is more
general than classical mathematics in the same way that noncommutative
algebra is more general than commutative algebra.  Secondly, even in
the setting of classical mathematics constructive reasoning can be
useful.  For example, it is possible to reason constructively in
Grothendieck toposes (which do not in general admit classical
reasoning).  Finally, proofs given in a constructive setting will
carry algorithmic content, whereas this is not true in general for
proofs given in the classical setting.\footnote{It should be mentioned
  that there are logical techniques, which themselves exploit the algorithmic content of constructive reasoning, for
  extracting algorithmic content from classical proofs.  See \cite{Kohlenbach:APT}.}

\subsection{Reasoning about spaces in type theory}

Voevodsky \cite{Vo2012a} has described a construction of \emph{set quotients}
of types.  Explicitly, a relation on a type $X$ is given by a map
$R:X\times X\to\hprop$.  For equivalence relations $R$, we can form
the quotient $X/R$ of $X$.  This type $X/R$ is necessarily a set and
can be shown to have the appropriate universal property expected by
such a quotient.  The set $\pi_{0}(X)$ of path components of $X$ is
constructed as a set quotient in the usual way.  Because loop spaces
$\Omega(X,x)$ of types $X$ at points $x:X$ can be defined type
theoretically, it then follows that we may construct all of the
higher-homotopy groups $\pi_{n}(X,x)$ of $X$ with basepoint $x:X$ by
setting 
\begin{align*}
  \pi_{n}(X,x) & := \pi_{0}\bigl(\Omega^{n}(X,x)\bigr).
\end{align*}
Many of the usual properties of the groups $\pi_{n}(X,x)$ can then be
verified type theoretically.  E.g.,
\begin{enumerate}
\item the homotopy groups of contractible spaces are $0$.
\item the usual Eckmann-Hilton \cite{Eckmann:1961tx} argument shows that $\pi_{n}(X,x)$ is abelian for $n>1$
  (Licata \cite{Licata} has given a proof of this in the proof assistant Agda).
\item Voevodsky has developed a large part of the theory of homotopy
  fiber sequences.  Using this it is possible to construct the long
  exact sequence associated to a fibration.
\end{enumerate}

The notion of inductive type described in Section \ref{sec:inductive_intro}. is the type
theoretic analogue of the notion of a free algebraic structure on
a signature (a list of generating operations together with their
arities) as studied in universal algebra.  By considering a type
theoretic analogue of free algebraic structures on a signature subject to
relations, it is possible to describe many familiar spaces type
theoretically.  This notion is that of \emph{higher-inductive
  type} which is currently being developed by a number of researchers
(cf., the work by Lumsdaine and Shulman \cite{Lumsdaine:HIT,Shulman:HTTVI}).
Rather than giving a comprehensive introduction to this
subject, we will give a simple example which should illustrate the
ideas and we will then summarize a few of the additional things which
can be done with this idea.

To describe the circle $S^{1}$ as a higher-inductive type, we require
that it should have one generator $b : S^{1}$ and one generator
$\ell:\id{S^{1}}(b,b)$.  One then obtains an induction principle for
$S^{1}$ similar to the induction principle for $\nat$ described
earlier.  Namely, given any type $\judge{x:S^{1}}{E(x)}$ fibered over
$S^{1}$ together with terms
\begin{itemize}
\item $b':E(b)$; and
\item $\ell':\id{E(b)}\bigl(\ell_{!}(b'),b'\bigr)$
\end{itemize}
there exists a term $\judge{x:S^{1}}{\rec_{S^{1}}(b',\ell',x):E(x)}$
(satisfying appropriate ``computation'' conditions).  (Here
$\ell_{!}(b')$ is $b'$ \emph{transported in the fiber along the loop $l$}
for which we refer the reader to Section \ref{sec:transport}.)  So, in
particular, in order to construct a map $S^{1}\to X$, it suffices to
give a point $x:X$ and a loop $\ell':\id{X}(x,x)$ on $x$.  The usual
properties of $S^{1}$ then follow from this type theoretic description
(see \cite{Shulman:S1Z} for a type theoretic proof that $\pi_{1}(S^{1},b)\iso\mathbb{Z}$).

In basically the same way, finite (and suitably inductively generated)
CW-complexes and relative CW-complexes can be constructed as
higher-inductive types.  In fact, Lumsdaine \cite{Lumsdaine:MSHIT} has shown that
with higher-inductive types, the syntax of type theory gives rise to
all of the structure of a model category except for the finite
limits and colimits.

\subsection{Future directions}

There are a number of exciting directions in homotopy type theory and
univalent foundations which are currently being pursued.  We will
briefly summarize several of them.

First, there are a number of interesting theoretical questions
surrounding the Univalence Axiom which remain open.  The most pressing
of these questions is the question of the ``constructivity'' of the
Univalence Axiom posed by Voevodsky in \cite{Vo2012}.  Voevodsky
conjectured that, in the presence of the Univalence Axiom, for every
$t:\nat$, there exists a numeral $n:\nat$ and a term of type
$\id{\nat}(t,n)$.  Moreover, it is expected that there exists an
algorithm which will return such data when given a term $t$ of type
$\nat$.  Additionally, there is the question of finding additional
models of the Univalence Axiom and characterizing such models (see,
e.g., \cite{Shulman:UAID}).

Secondly, it remains to be seen how much of modern homotopy theory can
be formalized in type theory using either higher-inductive types or
some other approach.  Indeed, there is still work to be done to arrive
at a complete theoretical understanding of higher-inductive types.

Finally, the formalization of ordinary (set level and higher-level)
mathematics in the univalent setting remains to be done.  At present,
a large amount of mathematics has been formalized by Voevodsky in his
Coq library.  Additionally, together with Voevodsky, the authors have
been working on developing an approach to the theory of integrable
systems (using the new notion of $p$-adic integrable system as a test
case) in the univalent setting \cite{PeVoWa2012}.  Ultimately, we hope that it
will be possible to formalize large amounts of modern mathematics in
the univalent setting and that doing so will give rise to both new
theoretical insights and good numerical algorithms (extracted from code in
a proof assistant like Coq) which can be applied to real world
problems by applied mathematicians.

\section{Basic Coq constructions} \label{sec2}

We will now introduce some basic constructions in Coq and their
corresponding homotopy theoretic interpretations.  We mention here
that there is an accompanying Coq file which includes all of the Coq
code discussed here, as well as some additional code.\footnote{The
  Coq file can be found either as supplementary data attached to the
  arXiv version of this paper or on the second author's webpage.}

\subsection{The Coq proof assistant} \label{Coq}

The Coq proof assistant \cite{Coq,Bertot:2004uj,Chlipala:CP} is a
computer system which is based one flavor of Martin-L\"of type theory
called the \emph{calculus of inductive constructions} and based in
part on the earlier \emph{calculus of constructions}
\cite{Coquand:1988eq} due to Coquand and Huet. 

In February 2010 Voevodsky \cite{Vo2012a} began writing a Coq library
of formalized mathematics based on the univalent model.  The resulting
library can currently be found online at the following location:
\begin{center}
\url{http://github.com/vladimirias/Foundations/}. 
\end{center}
There is also an HTML version of the library which can be found at Voevodsky's web page
\begin{center}
\url{http://www.math.ias.edu/~vladimir}
\end{center}
In addition to Voevodsky's library, there is also a repository which
is being developed by other researchers in homotopy type theory which
can be found at
\begin{center}
  \url{http://github.com/HoTT/HoTT}
\end{center}
Documentation on how to configure Coq for each of these libraries can
be found on the respective websites.  We expect a unification of these
libraries to occur at some point during the Special Year being held at
the Institute for Advanced Study during the 2012 - 2013 academic year.
For now though we work mostly in the style of Voevodsky's library.

While we will not explain here how to install or process a Coq file,
it is nonetheless worth mentioning that the way a Coq file is
generally processed is in an \emph{active} manner.  That is, one
processes the file in a step-by-step way and as one does so Coq
provides feedback regarding the current state of the file.
\begin{figure}[H]
\includegraphics[width=4.9in]{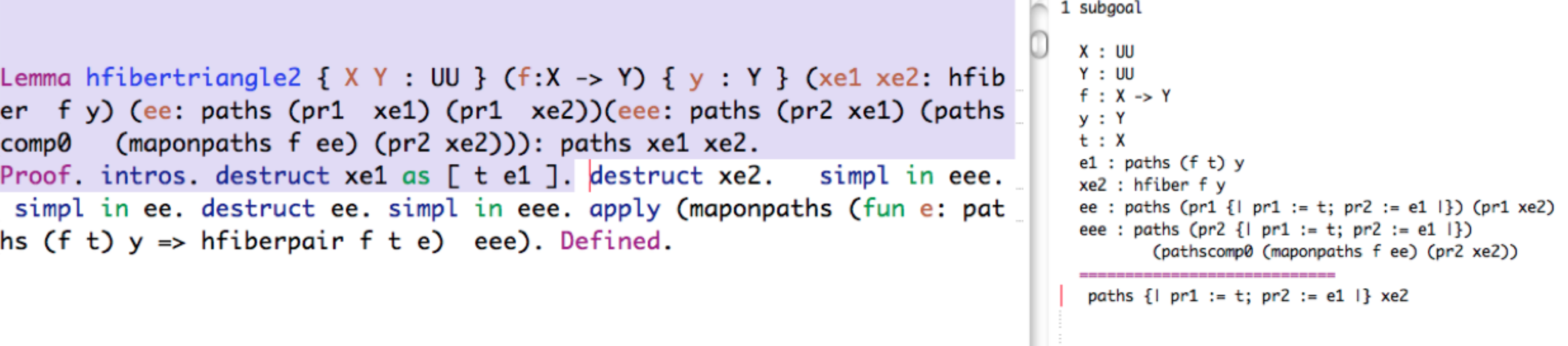}
  \caption{A screen with separate compartments when working with
    Coq. The code entered by the user is here displayed in the left
    compartment. The next goal required by Coq in order to complete
    the currently active proof appears in the upper right compartment.}
  \label{fig:diffeo}
\end{figure}
For example, Figure \ref{fig:diffeo} illustrates a Coq file currently
being processed.  The user is currently in the middle of processing a
proof (indicated in the left-hand pane of the image) and the shaded
text denotes the part of the file which has so far been processed by
Coq.  The right-hand pane is one of two mechanisms which Coq has for
providing the user with feedback.  In particular, this pane indicates
the current state of the proof which is being carried out.  Thus, as
the user progresses through a proof the output changes so as to always
indicate what remains to be done in order to complete the proof.
Further, more detailed, examples of this process are given below.

\subsection{Types and terms in Coq}\label{sec:sentences}

The Coq proof assistant, being based as it is on a form of type
theory, allows us to formalize and verify reasoning about types and
terms.  Throughout our discussion, the reader should have in mind the
interpretation of type theory described in Section
\ref{sec:homotopy_interp}.  Coq comes with a number of types and
type forming operations already built-in.  Using these it is possible
to define new types.  The first thing we want to do is to select a
fixed universe of small types with which we will work.  This is
accomplished by the following code:
\begin{center}
  \begin{coqcode}
Definition UU := Type.
  \end{coqcode}
\end{center}
Here the expression \verb|Type| is a built-in type in the Coq
system which is a universe of types (in a suitable technical sense).
The definition above then serves to define \verb|UU| to be this
fixed built-in universe of types.  We think of the terms of type
\verb|UU| as the small spaces and the universe \verb|UU|
itself as the (large) space of small spaces.
Mathematically this is roughly the same as fixing a Grothendieck
universe and letting \verb|UU| be the corresponding space of
spaces in the universe.  

The main reason for doing this, aside from notational convenience, is
a technical one arising from the internal
mechanism of Coq.  Namely, Coq secretly assigns indices to each
occurrence of \verb|Type| in a way which ensures a consistent
indexing.  However, we would like to work with one fixed universe and
not with an entire hierarchy thereof, and this is accomplished by
adopting the definition above.\footnote{If you would like to see the
  explicit indexing of universes, then you can add the line
  {\texttt{Set Printing Universes.}} to your Coq file.}  
The type \verb|UU| corresponds to $\UU$ from Section
\ref{sec:univalent_model} above.  Henceforth, any
statement of the form \verb|A : UU| should be thought of as
asserting that $A$ is a small space.

One interesting feature of the Coq system is that types are themselves
terms.  In particular, the ``type'' \verb|UU| above is itself a term of type \verb|Type|, where
this latter \verb|Type| is given the index $(n+1)$ when
\verb|UU| has index $n$. That is, being a type is really the same as being
a term in a higher universe.

\subsection{A direct definition involving function spaces}\label{sec:direct_def}

In order to illustrate some further features of the Coq system, we
will define some basic construction on function spaces.  First, we
define, for any small type, the identity function:
\begin{center}
  \begin{coqcode}
Definition idfun ( A : UU ) : A -> A := fun x => x.
  \end{coqcode}
\end{center}
Let us dissect this line of code and try to understand each of the
ingredients.  A definition, such as this one, is what we will call a
\emph{direct definition} and such a definition has the abstract form
summarized (together with the two examples we have so far encountered)
in Figure \ref{fig:direct}.
\begin{figure}[H]
  \centering
  \begin{tabular}{cccccccc}
    \verb|Definition| & \emph{name} &
    \emph{parameters} & : & \emph{type} & := &
    \emph{explicit definition} & .\\
    \hline
    & \verb|UU| & & & & & \verb|Type| &\\
    & \verb|idfun| & \verb|( A : UU )| &  &
    \verb| A -> A | & & \verb|fun x => x| &
  \end{tabular}
  \caption{Direct definitions in Coq.}
  \label{fig:direct}
\end{figure}

Several remarks about Figure \ref{fig:direct} are in order.  First, the \emph{name} is the name
given to the term.  This can be whatever (modulo some restrictions on
the syntactic form) the user likes.  The \emph{type} is the type of
the term being defined.  I.e., we have that \emph{name} is of type
\emph{type}.  The next thing to note is that the \emph{parameters} can
be a list of terms variables of fixed types.  In the case of
\verb|idfun| there is just a single parameter: the type
\verb|A : UU|; in the case of \verb|UU| there are no
parameters at all.  Within a definition, the parameters should be
enclosed in brackets as in \verb|( A : UU )|.  Next, note that
it is not strictly necessary to declare the type.  When no type is
given, Coq will infer the type.  Finally, note that the period at the
end of the definition must be included in order for Coq to correctly
parse the input.

Coming back to the definition of \verb|idfun|, it is worth
mentioning that the type \verb|A -> A| is the way of denoting the
function space $A^{A}$ in Coq.  That is, for types \verb|A| and
\verb|B|, the type \verb|A -> B| is the type of functions
from \verb|A| to \verb|B|.  For us, this type should be
thought of more specifically as the type of all continuous functions
from the space \verb|A| to the space \verb|B|.  The
remaining part of this definition is the actual content of the
definition: \verb|fun x => x|.  In this definition, the
expressions \verb|fun| and \verb|=>| go together and tell us
that it is the function which takes a point $x$ in $A$ and gives back $x$.
That is, \verb|fun ... => ...| is the same as giving a definition
of a function by writing $\ldots\mapsto\ldots$ or (for those familiar
with lambda calculus) using lambda abstraction.  In particular, the
definition states that \verb|idfun| is the function given by
$x\mapsto x$ (i.e., $\lambda_{x}.x$).

We can now play around a bit with the type checking mechanisms of
Coq.  Let us enter the following into the Coq code:
\begin{center}
  \begin{coqcode}
Section idfun_test.    
Variable A : UU.
Variable a : A.
  \end{coqcode}
\end{center}
The first line tells Coq that we are starting a new section of the
file in which we will introduce certain hypotheses.  The next line
tells Coq that we would like to assume, for the duration of the
current section, that \verb|A| is a space in \verb|UU|.  The
final line similarly tells Coq that we are assuming given a point
\verb|a| of \verb|A|.  Now, if we add the following line to
our Coq file and process up to this point, we will be able to see what
type the term \verb|idfun A| has:
\begin{center}
  \begin{coqcode}
Check idfun A.
  \end{coqcode}
\end{center}
Coq will respond by telling us 
\begin{center}
  \begin{coqcode}
idfun A : A -> A.
  \end{coqcode}
\end{center}
Similarly, if we enter
\begin{center}
  \begin{coqcode}
Check idfun _ a.
  \end{coqcode}
\end{center}
Then Coq replies with \verb|idfun A a : A|.  Here note that we
write \verb|idfun _ a| to tell Coq that we would like for it to
guess the parameter (in this case \verb|A|) which should go
in the place indicated by the underscore.  Finally, we close the
section by entering
\begin{center}
  \begin{coqcode}
End idfun_test.
  \end{coqcode}
\end{center}
After entering this line of code, the variables \verb|A| and
\verb|a| are no longer declared.\footnote{You can verify this by
  trying \texttt{Check A.} and observing the response from Coq.
  Be sure to remove this line from your code though otherwise you
  won't be able to go any further!}

\subsection{An indirect definition involving function spaces}

We will now show, given functions $f:A\to B$ and $g:B\to C$, how to
construct the composite $g\circ f : A \to C$ type theoretically.  In
order to introduce \emph{indirect} definitions, we will give two
ways to construct $g\circ f$.

The utility of indirect definitions in Coq is that sometimes it is not
easy to see how to give the explicit definition of a term.  This is
especially true as one starts working with increasingly complicated
definitions.  As such, rather than having to struggle to define
exactly the required term it is possible to construct the term being
defined as a kind of proof.  Along the way, as this proof is
constructed, certain automation possible in Coq can be employed.  In
order to see how this works in practice, let us introduce our first
indirect definition.
\begin{center}
  \begin{coqcode}[frame=single,backgroundcolor=\color{mylight},framerule=0pt]
Definition funcomp_indirect ( A B C : UU ) ( f : A -> B ) ( g : B -> C ) : A -> C.
Proof.
  intros x. apply g. apply f. assumption.
Defined.
  \end{coqcode}
\end{center}
The first observation about this definition is that it looks like
everything to the left of the \verb|:=| in a direct definition.
In this case there are three parameters of type \verb|UU|
(namely, \verb|A|, \verb|B| and \verb|C|).  There is
also one parameter of type \verb|A -> B| (namely,
\verb|f|).  Finally, there is one parameter of type \verb|B -> C| (namely, \verb|g|).  

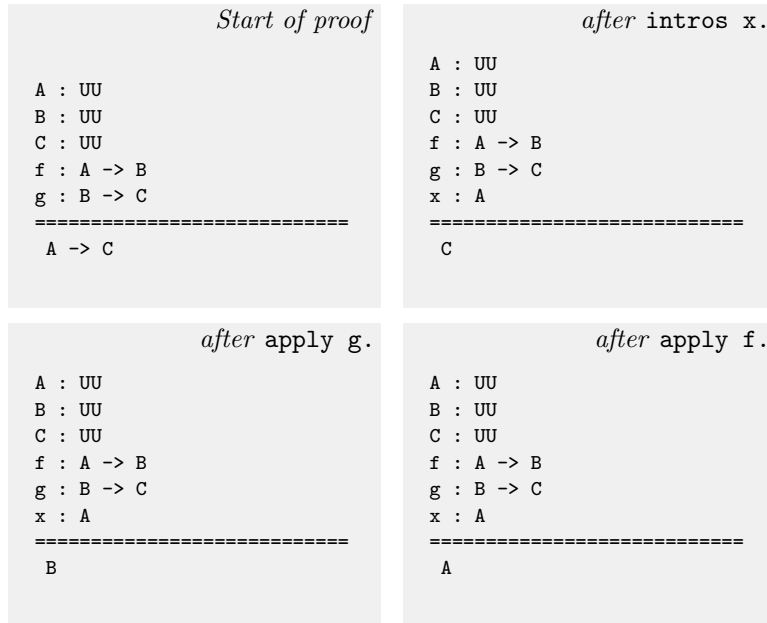
\begin{figure}[H]
  \begin{tikzpicture}
    \node[smallcoqbox] (zero)  at (0,0) {%
      \begin{minipage}{4.25cm}
        \footnotesize
        
        \vphantom{\texttt{x : A}}
        
        \noindent\verb|A : UU|
        
        \noindent\verb|B : UU|
        
        \noindent\verb|C : UU|
        
        \noindent\verb|f : A -> B|
        
        \noindent\verb|g : B -> C|
        
        \noindent\verb|============================|
        
        \noindent\verb| A -> C|
      \end{minipage}
    };
    \node[anchor=north east, inner sep=2pt] (titlezero) at
    (zero.north east) {\emph{Start of proof}};
    \node[smallcoqbox] (one) at (5.25,0) {%
      \begin{minipage}{4.25cm}
        \footnotesize
        \noindent\verb|A : UU|
        
        \noindent\verb|B : UU|
        
        \noindent\verb|C : UU|
        
        \noindent\verb|f : A -> B|
        
        \noindent\verb|g : B -> C|
        
        \noindent\verb|x : A|
        
        \noindent\verb|============================|
        
        \noindent\verb| C|
      \end{minipage}
    };
    \node[anchor=north east, inner sep=2pt] (titleone) at
    (one.north east) {\emph{after} \verb|intros x.|};
    \node[smallcoqbox] (two) at (0,-4.25) {%
      \begin{minipage}{4.25cm}
        \footnotesize
        \noindent\verb|A : UU|
        
        \noindent\verb|B : UU|
        
        \noindent\verb|C : UU|
        
        \noindent\verb|f : A -> B|
        
        \noindent\verb|g : B -> C|
        
        \noindent\verb|x : A|
        
        \noindent\verb|============================|
        
        \noindent\verb| B|
      \end{minipage}
    };
    \node[anchor=north east, inner sep=2pt] (titletwo) at
    (two.north east) {\emph{after} \verb|apply g.|};
    \node[smallcoqbox] (three) at (5.25,-4.25) {%
      \begin{minipage}{4.25cm}
        \footnotesize
        \noindent\verb|A : UU|
        
        \noindent\verb|B : UU|
        
        \noindent\verb|C : UU|
        
        \noindent\verb|f : A -> B|
        
        \noindent\verb|g : B -> C|
        
        \noindent\verb|x : A|
        
        \noindent\verb|============================|
        
        \noindent\verb| A|
      \end{minipage}
    };
    \node[anchor=north east, inner sep=2pt] (titlethree) at
    (three.north east) {\emph{after} \verb|apply f.|};
  \end{tikzpicture}
  \caption{Coq output during indirect definition of function
    composition.}
  \label{figure:indirect_funcomp}
\end{figure}

After the first full stop of an indirect definition, we encounter the
start of the proof.  This is given by the line
\begin{center}
  \begin{coqcode}
Proof.
  \end{coqcode}
\end{center}
Likewise, the end of the proof is indicated by 
\begin{center}
  \begin{coqcode}
Defined.
  \end{coqcode}
\end{center}
Between the start of the proof and the end of the proof is a sequence
of what are called \emph{tactics}, which allow one to construct, using
the given parameters, the required term.  One limitation of writing an
article which includes proofs in Coq, is that proofs in Coq are
usually constructed using ``backward'' reasoning and so it can be hard to
read for the uninitiated.  In particular, the nature of Coq is such
that, \emph{qua} interactive proof assistant, proofs can be understood
better by directly watching the output of a Coq session, where an
additional window appears after each step, giving us precise
explanations on any given step of the proof.  We have included in
Figure \ref{figure:indirect_funcomp} the output from Coq as we move
through the proof of \verb|funcomp_indirect|.  Readers
should not be discouraged if they are unable to read Coq proofs
directly: it is much easier if you are going through the proofs yourself in the
computer.

We will now go through the proof one step at a time.  Notice that once the proof is started Coq
displays the hypotheses (above the line \verb|=======|) together
with the current goal (below the line).  Since the goal is to
construct a function \verb|A -> C| we are allowed to assume given an
arbitrary term \verb|x : A|.  This is accomplished in Coq by
entering \verb|intros x|.  Note that the name \verb|x| here
is something which we have chosen and the user can choose this name
freely (or it can be omitted, in which case the Coq system will supply
a name of its own choosing).  As such, after processing
\verb|intros x|, the output has changed (as indicated in Figure
\ref{figure:indirect_funcomp}) by adding a new hypothesis 
(\verb|x : A|) and changing the goal to \verb|C|.  This means that we now
need to supply a term of type \verb|C|.  To accomplish this, we
note that we are given a function \verb|g : B -> C| and so it
suffices to supply a term of type \verb|B|.  We communicate the
fact that we intend to use \verb|g| to obtain the goal by
entering \verb|apply g|.  The effect of this is to change the
goal from \verb|C| to \verb|B| since \verb|B| is the
domain of \verb|g|.  Applying the same reasoning now with
\verb|f| we are in the final situation indicated in Figure
\ref{figure:indirect_funcomp}.  Because we have as a hypothesis the term
\verb|x : A| and the current goal is to construct a term of type
\verb|A| we may simply communicate to the Coq system that
there is already a term of the required type appearing among the
hypotheses.  This is accomplished by entering \verb|assumption|.
Indeed, at this point, Coq tells us
\begin{center}
  \begin{coqcode}
No more subgoals.
  \end{coqcode}
\end{center}
and the proof is complete.  Note that we \emph{must} add the final
``\verb|Defined.|'' in order for the Coq system to correctly
record the proof.

We now turn to the corresponding direct definition:
\begin{center}
  \begin{coqcode}
Definition funcomp { A B C : UU } ( f : A -> B ) ( g : B -> C ) := fun x : A => g ( f x ).
  \end{coqcode}
\end{center}
One thing worth noting about this definition is that here we have
enclosed the first three parameters in curved brackets as \verb|{A B C : UU}| in order to indicate to the Coq system that these
parameters are \emph{implicit}.  Implicit parameters do not need to be
supplied (when the term is applied) by the user and the system will
try to infer the values of these parameters.  In this case, these can
be inferred from the types of \verb|f| and \verb|g|.  Note
also that we have here not given explicitly the type of the term being
defined.  As such, we must include the additional typing data
\verb|x : A| in the definition in order for the Coq system to be
able to infer the type of the term. 

We can now check that our definition agrees with the direct
one by entering:
\begin{center}
  \begin{coqcode}
Print funcomp.
Print funcomp_indirect.
  \end{coqcode}
\end{center}
The effect of \verb|Print| is to output both the type and
explicit definition of the term in question.  In particular, even if
the term in question was defined indirectly, as our
\verb|funcomp_indirect| was, it is an \emph{explicit term
as far as Coq is concerned}, and when \verb|Print| is used Coq
will unfold the term to give a completely explicit description.

To understand what this means, simply think of a linear differential
equation for which it is possible to explicitly write down
the solutions. The solution set can be difficult to write down, but it can
be done. Although you may not want to have these solutions in front of
you, you know that the explicit solutions are available if they are required (say, to
verify that they satisfy a certain equation).  The Coq system
effectively keeps track of this kind of book-keeping for the user.

\section{Some basic inductive types} 
\label{sec:inductive}

First we explain two important notions: induction
and recursion, and how they relate to each other. As a warm up,
we are going to explain them for the case of the natural numbers,
but it is important to keep in mind that for us, the key inductive
and recursive process will take place over more general objects
than the natural numbers. But to understand those, first one
must understand the simpler case of the natural numbers.

\subsection{The inductive type of natural numbers}\label{sec:induction}

Although it is possible to construct arbitrary inductive types in Coq,
we will start by looking at one inductive type which is already
defined in the Coq system.  This is the type \verb|nat| of
natural numbers.  As an inductive type, \verb|nat| is generated
by 
\begin{itemize}
\item the single generator ($0$-ary operation) \verb|0 : nat|; and
\item the single generating function ($1$-ary operation) 
  \verb| S : nat -> nat | (this is the usual \emph{successor function}).
\end{itemize}
In Coq, this inductive type is specified as follows:
\begin{center}
  \begin{coqcode}
Inductive nat := 0 | S : nat -> nat.
  \end{coqcode}
\end{center}
Here the internal Coq command \verb|Inductive| functions
similarly to \verb|Definition| (as we will see below).  For now
the crucial point is to observe that the generating operations of the
inductive type appear on the right of \verb|:=| and are separated
by the symbol $|$.

One of the advantages of working with an inductive type such as the
natural numbers is that functions with inductive domain can be defined
by cases.  To see an example, consider the predecessor function:
\begin{center}
\begin{coqcode}
Definition predecessor ( n : nat ) : nat := 
 match n with
   | 0 => 0
   | S m => m
 end.
\end{coqcode}
\end{center}
This is the way of telling Coq that the predecessor function is the
function $\text{predecessor}:\nat\to\nat$ given by case analysis as
\begin{align*}
  \text{predecessor}(n) & :=
   \begin{cases}
    0 & \text{ if }n=0,\text{ and}\\
    m & \text{ if }n=(m+1).
  \end{cases}
\end{align*}
Because we are using the type \verb|nat| which is predefined in
the Coq system, Coq will recognize that ordinary numerals refer to the
corresponding terms of type \verb|nat|.  So, for example, Coq
knows that \verb|3| is the same as \verb|S ( S ( S 0 ) )|.
Using this, we can play around with the computational abilities of
Coq by entering, say,
\begin{center}
  \begin{coqcode}
Eval compute in predecessor 26.
  \end{coqcode}
\end{center}
Here \verb|Eval compute in| tells Coq that we would like it to
compute the value of the subsequent expression (in this case
\verb|predecessor 26|).  Coq correctly replies
\begin{center}
  \begin{coqcode}
= 25 : nat
  \end{coqcode}
\end{center}
as expected.

\begin{figure}[H]
  \begin{tikzpicture}
    \node[smallcoqbox] (zero)  at (0,0) {%
      \begin{minipage}{4.25cm}
        \footnotesize
        \noindent\verb|n : nat|
        
        \noindent\verb|============================|
        
        \noindent\verb| nat|
      \end{minipage}
    };
    \node[anchor=north east, inner sep=2pt] (titlezero) at
    (zero.north east) {\emph{Start of proof}};
    \node[smallcoqbox] (one) at (5.25,0) {%
      \begin{minipage}{4.25cm}
        \footnotesize
        \vphantom{\texttt|n : nat |}
        
        \noindent\verb|============================|
        
        \noindent\verb| nat|
      \end{minipage}
    };
    \node[anchor=north east, inner sep=2pt] (titleone) at
    (one.north east) {\emph{after} \verb|destruct n.|};
    \node[smallcoqbox] (two) at (0,-2.5) {%
      \begin{minipage}{4.25cm}
        \footnotesize
        \noindent\verb|n : nat|
        
        \noindent\verb|============================|
        
        \noindent\verb| nat|
      \end{minipage}
    };
    \node[anchor=north east, inner sep=2pt] (titletwo) at
    (two.north east) {\emph{after} \verb|exact 0.|};
  \end{tikzpicture}
  \label{figure:indirect_predecessor}
  \caption{Coq output during indirect definition of the predecessor function.}
\end{figure}
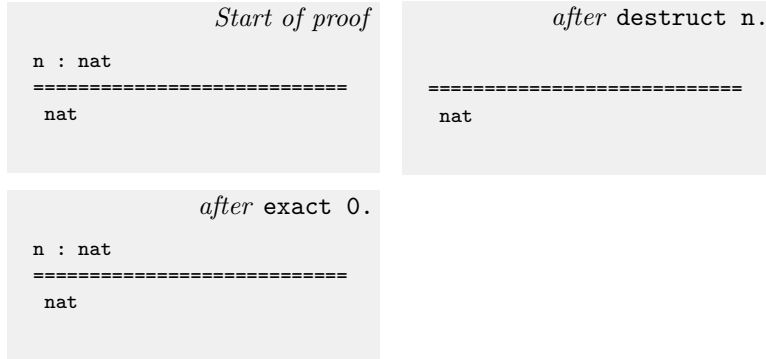

The definition of the predecessor function given above using
\verb|match| is a direct definition as described above in
Section \ref{sec:direct_def}.  It is also possible to define the
predecessor function via an indirect definition as follows:
\begin{center}
  \begin{coqcode}
Definition indirect_predecessor ( n : nat ) : nat.
Proof.
  destruct n. exact 0. exact n.
Defined.
  \end{coqcode}
\end{center}
In the proof there are two new tactics.  The first,
\verb|destruct n|, tells the Coq system that we will reason by
cases on the structure of \verb|n| as a term of type
\verb|nat|.  Coq knows that, as a natural number, there are two
cases and in the first case there is no hypothesis necessary (see
Figure \ref{figure:indirect_predecessor}) because \verb|n| is
\verb|0| in this case.  At this stage, we know that we would like
the output of the function to be \verb|0| and we tell Coq this
using the \verb|exact| tactic.  In general, if we enter
\verb|exact| $x$, this tells Coq that the term we are looking for
is \emph{exactly} the term $x$.  Once we have entered
 \verb|exact 0|, Coq moves on to the second possibility: the term is a successor.
Note that in the list of hypotheses at this stage (see Figure
\ref{figure:indirect_predecessor}) the term \verb|n : nat| is
listed and so \emph{superficially} things are just as they were at the
start of the proof.  However, because we earlier employed the
\verb|destruct| tactic, Coq knows that we must now give the
required output of the predecessor function when given the value
\verb|S n|.  As such, we enter \verb|exact n|.  Comparing
\verb|predecessor| and \verb|indirect_predecessor| using the
\verb|Print| command reveals that they are indeed identical terms.
We will now turn to several of the inductive types more closely
related to the homotopy theoretic side of things.

\subsection{Fibrations and the total space of a fibration} 

Fibrations can be understood as a homotopy theoretic generalization of
the more familiar notion of fiber bundle.  Similarly, they can be
understood as a homotopy theoretic version of Grothendieck fibrations
familiar from category theory and algebraic geometry.  In particular,
for spaces, a fibration is a map
\begin{align*}
\pi \colon E \to B
\end{align*}
which possesses a certain homotopy-lifting property.  In this case we
would refer to $B$ as the \emph{base space} and to $E$ as the
\emph{total space} of the fibration.  Given a point $b$ in $B$, the
\emph{fiber} over $b$, which we sometimes write as $E_{b}$, is just the preimage $\pi^{-1}(b)$ of $b$ under
the map $\pi$.  As mentioned in Section \ref{sec:homotopy_interp}
above, fibrations correspond to
types which depend on parameters (so-called \emph{dependent types}).
In Coq, the way of dealing with dependent types is somewhat different
from the one sketched in Section \ref{sec:homotopy_interp}.  In
particular, given an element \verb|B : UU| a fibration for us is
a term
\begin{center}
  \begin{coqcode}
E : B -> UU.
  \end{coqcode}
\end{center}
The idea that a fibration over a base space $B$ is the same as a map
from $B$ into a suitable universe is a classical idea, especially in
category theory where it is a basic fact that Grothendieck fibrations
over a category $B$ correspond to pseudo-functors from $B$ into the
2-category of small categories.  (These ideas are ubiquitous in
category theory and have been developed in great generality by the
Australian school of category theorists.  A nice exposition of these
ideas can be found in the first section of \cite{Street:1980da}.)

The idea behind this correspondence is that a fibration can be completely recovered from its base space
$B$ together with its fibers by gluing the fibers together in a
coherent way in accordance with the structure of the base space.  The
same can be accomplished in Coq by defining the total space of a
fibration \verb|E : B -> UU| as an inductive type.
\begin{center}
  \begin{coqcode}
Inductive total {B :UU} (E : B -> UU) : UU := 
pair ( x : B ) ( y : E x ).
  \end{coqcode}
\end{center}
Intuitively, \verb|total E| should be thought of as a space consisting of all
pairs \verb|(b,e)|, where \verb|b| is a point of the base
space \verb|B| and \verb|e| in a point of the fiber
\verb|E b| over \verb|b|.

Fiber bundles $E\to B$ are sometimes thought of as a ``twisted'' generalization
of direct products $F\times B\to B$, and the fact that fibrations are
a homotopical generalization of this notion reveals itself type
theoretically by the fact that the total space
construction \verb|total| is a generalization of the construction
of direct products of spaces $A\times B$, which are given by
\begin{center}
  \begin{coqcode}
Definition dirprod { A B : UU } : UU := 
 total ( fun x : A => B ).
  \end{coqcode}
\end{center}
Returning to \verb|total|, we define a projection map
\verb|total E -> B| by 
\begin{center}
  \begin{coqcode}
Definition pr1 { B : UU } { E : B -> UU } : total E -> B := fun z => match z with pair x y => x end.
  \end{coqcode}
\end{center}
This map serves to exhibit \verb|E| as a fibration over
\verb|B|.

\section{The path space} \label{sec:pathspace}

As mentioned in Sections \ref{sec:inductive_intro} and
\ref{sec:groupoids}, the type theoretic \emph{identity type} is
interpreted homotopy theoretically as the path space.  The path space,
and the corresponding type in Coq, is so important that we will now
carefully describe several basic constructions involving it in the
setting of Coq.
\begin{figure}[H]
  \centering
  \begin{tikzpicture}[color=mydark, fill=mylight, line width=1pt]
    \draw[fill=mylight] plot [smooth cycle,tension=.75] coordinates
    {(-2.5,0) (-2,1.5) (0,1.25) (1.75,1.5) (2.75,0) (1.75,-1) (-1,-1)};
    \draw[fill=white] plot [smooth cycle,tension=.75] coordinates
    {(0,.2) (.5,.45) (1,.2) (.5,-.05) };
    \node[circle,fill=mydark,inner sep=0pt,outer sep=4pt,minimum size=.75mm] (aa)
    at (.2,-.4) {};
    \node at (.15,-.6) {$\scriptstyle a$};
    \node[circle,fill=mydark,inner sep=0pt,outer sep=4pt,minimum size=.75mm] (bb)
    at (-1,.75) {};
    \node at (-1.2,.9) {$\scriptstyle b$};
    \draw[color=mydark,->-=.5,line width=.5pt] plot [smooth,tension=.75]
    coordinates { (aa) (1.2,.2) (.4,.8) (bb) };
    \node[font=\tiny] at (1.3,.2) {$k$};
    \setcounter{myi}{0}
    \foreach \i in {1,...,9}
    {
      \pgfmathsetcounter{myi}{\themyi+9.9}
      \setcounter{myi}{\themyi}
      \draw [mydark,line width=.1pt] plot [smooth,tension=.5] 
      coordinates { (aa) (0-.1*\i,0+.01*\i) (bb) }; 
    }
    \draw [mydark,line width=.5pt,->-=.7] plot [smooth,tension=.5] 
    coordinates { (aa) (-1,.1) (bb) };
    \draw[mydark,line width=.5pt,->-=.4] plot [smooth,tension=.5] coordinates { (aa) (-.1,.01) (bb) };
    \node[circle,inner sep=.2mm,fill=mylight,fill opacity=.4,text opacity=1] at (-.7,.2) {$\scriptstyle p$};
    \node[font=\tiny] at (-.25,.4) {$f$};
    \node[font=\tiny] at (-.8,-.2) {$g$};
    \draw[dotted,line width=.6pt] (-1,3) to (bb);
    \draw[fill=mylight!50,line width=.8pt] plot [smooth cycle,tension=-.15] coordinates
    {(-2.2,2.5) (-2.2,4.5) (.2,4.5) (.2,2.5) };
    \node[circle,fill=mydark,inner sep=0pt,outer sep=2pt,minimum
    size=.75mm] (ee) at (-.5,3.4) {};
    \node[circle,fill=mydark,inner sep=0pt,outer sep=2pt,minimum
    size=.75mm] (ff) at (-1.2,3.9) {};
    \node[circle,fill=mydark,inner sep=0pt,outer sep=2pt,minimum
    size=.75mm] (gg) at (-1.4,3.1) {};
    \draw[mydark,line width=.4pt,->-=.5,inner sep=0pt, outer sep=0pt] (-1.2,3.9) to (-1.4,3.1);
    \node[font=\tiny] at (-.38,3.5) {$\scriptstyle k$};
    \node[font=\tiny] at (-1.05,4) {$\scriptstyle f$};
    \node[font=\tiny] at (-1.5,3) {$\scriptstyle g$};
    \node[font=\tiny] at (-1.45,3.5) {$\scriptstyle p$};
    \node[font=\tiny] at (.8,4.3) {$\texttt{paths }a\;b$};
    \node[font=\footnotesize] at (-2.55,1.45) {$A$};
  \end{tikzpicture} \label{fig:paths}
  \caption{The path space fibration $\texttt{paths }a$ with the fiber
    over a point $b$.  Here $p$ is a (path) homotopy from $f$ to $g$.}
\end{figure}
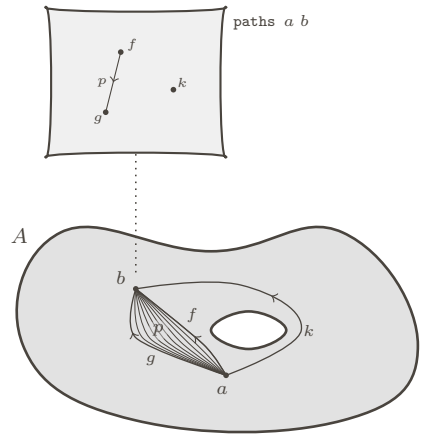
In Coq, the path space \verb|paths| is defined as follows.
\begin{center}
\begin{coqcode}
Notation paths := identity.
\end{coqcode}
\end{center}
Here \verb|identity|, like \verb|nat|, is a built-in
inductive type in the Coq system.  We can see how it is defined
inductively using \verb|Print| to find
\begin{center}
  \begin{coqcode}
Inductive identity (A : Type) (a : A) : A -> Type :=
 identity_refl : identity a a.
  \end{coqcode}
\end{center}
That is, for each \verb|a : A|, \verb|identity a| is the
fibration freely generated by a term \verb|identity_refl a| in
the fiber over \verb|a|.  

We add the following line in order to introduce a
slightly shorter notation for the terms \verb|identity_refl|:
\begin{center}
  \begin{coqcode}
Notation idpath := identity_refl.
  \end{coqcode}
\end{center}
That is, for \verb|a : A|, \verb|idpath a : paths a a| is
the \emph{identity path} based at \verb|a|.

Recall that a \emph{path} in a space $A$ is a continuous function
$\varphi:I\to A$ where $I=[0,1]$ is the unit interval.  We say that
$\varphi$ \emph{is a path from a point $a$ of $A$ to a point $b$ of $A$} when
$\varphi(0)=a$ and $\varphi(1)=b$.  Then, the \emph{path space} $A^{I}$ is
the space of paths in $A$ and it comes equipped with two maps
$\partial_{0},\partial_{1}:A^{I}\to A$ given by
$\partial_{i}(\varphi):=\varphi(i)$ for $i=0,1$.  The induced map
$\langle\partial_{0},\partial_{1}\rangle: A^{I}\to A\times A$ is a
fibration which gives a factorization of the diagonal $\Delta:A\to
A\times A$ as 
\begin{center}
  \begin{tikzpicture}[auto]
    \node (UL) at (0,1.25) {$A$};
    \node (UR) at (2.5,1.25) {$A^{I}$};
    \node (B) at (1.25,0) {$A\times A$};
    \draw[->] (UL) to (UR);
    \draw[->,bend right=10pt] (UL) to node[mylabel,swap] {$\Delta$} (B);
    \draw[->,bend left=10pt] (UR) to (B);
  \end{tikzpicture}
\end{center}
where the first map $A\to A^{I}$ is a weak equivalence and the second is the
fibration mentioned above.  Here the first map $A\to A^{I}$ sends a
point $a$ to the constant loop based at $a$.  (That is, this first map
is precisely \verb|idpath|.)  One of the many important contributions of Quillen in
\cite{Quillen:1967uz} was to demonstrate that it is in fact possible
to do homotopy theory without the unit interval provided that one has
the structure of path spaces, weak equivalences, fibrations, and a few
other ingredients.  This is part of the reason that, even though type
theory does not (without adding higher-inductive types or something
similar) provide us with a unit interval, it is still possible to work
with homotopy theoretic structures type theoretically.

\subsection{Groupoid structure of the path space}

We will now describe the groupoid structure which the path space
endows on each type.  These constructions are well-known and their
connection with higher-dimensional groupoids was first noticed by
Hofmann and Streicher \cite{Hofmann:1998ty}.

First, given a path $f$ from $a$ to $b$ in $A$ we would like to be
able to reverse this path to obtain a path from $b$ to $a$.  For
topological spaces this is easy because a path $\varphi:I\to A$ gives
rise to an inverse path $\varphi'$ given by
$\varphi'(t):=\varphi(1-t)$, for $0\leq t\leq 1$.
\begin{center}
  \begin{coqcode}
Definition pathsinv { A : UU } { a b : A } ( f : paths a b ) 
: paths b a.
Proof.
  destruct f. apply idpath. 
Defined.
  \end{coqcode}
\end{center}
Here recall that \verb|destruct| allows us to argue by cases
about terms of inductive types.  Here \verb|f| is of type
\verb|paths a b|, which is inductive, and therefore this tactic
applies.  In this case, there is only one case to consider:
\verb|f| is really the identity path 
\verb|idpath a : paths a a|.  Because the inverse of the identity is the identity we then use
\verb|apply idpath| to complete the proof.  (Note that we could
also have used \verb|exact ( idpath a )| instead of
\verb|apply idpath| here to obtain the same term.)

Next, given a path $f$ as above together with another path $g$ from
$b$ to $c$, we would like to define the composite path from $a$ to $c$
obtained by first traveling along $f$ and then traveling along $g$.
This operation of \emph{path composition} is defined as follows:
\begin{center}
  \begin{coqcode}
Definition pathscomp { A : UU } { a b c : A } ( f : paths a b ) ( g : paths b c ) : paths a c.
Proof.
  destruct f. assumption.
Defined.
  \end{coqcode}
\end{center}
Once again, the proof begins with \verb|destruct f| which
effectively collapses \verb|f| to a constant loop.  In
particular, the result of this is to change the ambient hypotheses so
that \verb|g| is now of type \verb|paths a c| (see Figure
\ref{figure:pathscomp}).  At this stage, the goal matches the type of
\verb|g| and we use \verb|assumption| to let the Coq system
choose \verb|g| as the result of composing \verb|g| with the
identity path.
\begin{figure}[H]
  \begin{tikzpicture}
    \node[smallcoqbox] (zero)  at (0,0) {%
      \begin{minipage}{4.25cm}
        \footnotesize
        \noindent\verb|A : UU|

        \noindent\verb|a : A|

        \noindent\verb|b : A|
        
        \noindent\verb|c : A|
        
        \noindent\verb|f : paths a b|
        
        \noindent\verb|g : paths b c|

        \noindent\verb|============================|

        \noindent\verb| paths a c|
      \end{minipage}
    };
    \node[anchor=north east, inner sep=2pt] (titlezero) at
    (zero.north east) {\emph{Start of proof}};
    \node[smallcoqbox] (one) at (5.25,0) {%
      \begin{minipage}{4.25cm}
        \footnotesize
        \vphantom{\texttt{b : A}}
        
        \vphantom{\texttt{f : paths a b}}

        \noindent\verb|A : UU|

        \noindent\verb|a : A|
        
        \noindent\verb|c : A|
        
        \noindent\verb|g : paths a c|

        \noindent\verb|============================|

        \noindent\verb| paths a c|
      \end{minipage}
    };
    \node[anchor=north east, inner sep=2pt] (titleone) at
    (one.north east) {\emph{after} \verb|destruct f.|};
  \end{tikzpicture}
  \caption{Coq output during the definition of path composition.}
  \label{figure:pathscomp}
\end{figure}
One immediate consequence of this definition is that the unit law
$f\circ 1_{a}=f$ for $f:a\to b$ holds \emph{on the nose} in the sense
that these two terms (\verb|pathscomp ( idpath a ) f| and 
\verb|f|) are identical in the strong $=$ sense.  On the other
hand, the unit law $1_{b}\circ f=f$ does not hold on the nose.
Instead, it only holds up to the existence of a higher-dimensional
path as described in the following Lemma:
\begin{center}
  \begin{coqcode}
Lemma isrunitalpathscomp { A : UU } { a b : A } ( f : paths a b ) : paths ( pathscomp f ( idpath b ) ) f.
Proof.
  destruct f. apply idpath. 
Defined.
  \end{coqcode}
\end{center}
The proof of this requires little comment (when $f$ becomes itself an
identity path, the composite becomes, by the left-unit law mentioned
above, an identity path).  The one thing to note here is that here
instead of \verb|Definition| we have written \verb|Lemma|.
Although there are some technical differences between these two ways
of defining terms they are for us entirely interchangeable and
therefore we use the appellation ``Lemma'' in keeping with the
traditional mathematical distinction between definitions and lemmas.

That facts that, up to the existence of higher-dimensional paths,
composition of paths is associative and that the inverses given
by \verb|pathsinv| are inverses for composition are recorded as
the terms \verb|isassocpathscomp|, \verb|islinvpathsinv| and
\verb|isrinvpathsinv|.  However, the descriptions of
these terms are omitted in light of the fact that they all
follow the same pattern as the proof of \verb|isrunitalpathscomp|.

\subsection{The functorial action of a continuous map on a path}

Classically, given a continuous map $f:A\to B$ and a path
$\varphi:I\to A$ in $A$, we obtain a corresponding path in $B$ by
composition of continuous functions.  Thinking of spaces as
$\infty$-groupoids, this operation of going from the path $\varphi$ in
$A$ to the path $f\circ\varphi$ in $B$ is the functorial action of $f$
on $1$-cells of the $\infty$-groupoid $A$.  In Coq, this action of
transporting a path in $A$ to a path in $B$ along a continuous map is
given as follows:
\begin{center}
  \begin{coqcode}
Definition maponpaths { A B : UU } ( f : A -> B ) { a a' : A } ( p : paths a a' ) : paths ( f a ) ( f a' ).
Proof.
  destruct p. apply idpath. 
Defined.
  \end{coqcode}
\end{center}
The proof again follows the familiar pattern: when the path $p$ is the
identity path on $a$, the result of applying $f$ should be the
identity path on $f(a)$.  We introduce the following notation for
\verb|maponpaths|:
\begin{center}
  \begin{coqcode}
Notation "f ` p" := ( maponpaths f p ) (at level 30 ).
  \end{coqcode}
\end{center}
This is an example of a general mechanism in Coq for defining
notations, but discussion of this mechanism is outside of the scope of
this article (the crucial point here is that the value 30 tells how
tightly the operation $`$ should bind).
\begin{figure}[H]
  \centering
  \begin{tikzpicture}
    \draw[fill=mylight] plot [smooth cycle,tension=.75] coordinates
    {(-4,0) (-3.5,1.5) (-2,1) (-.5,1.5) (0,0) (-2,-1) };
    \draw[fill=mylight] plot [smooth cycle,tension=1] coordinates
    {(2.5,0) (4,1) (5.5,0) (4,-1) };
    \node[circle,fill=mydark,inner sep=0pt,outer sep=4pt,minimum
    size=.75mm] (aa) at (-3.5,.5) {};
    \node[circle,fill=mydark,inner sep=0pt,outer sep=4pt,minimum
    size=.75mm] (aaprime) at (-.5,.5) {};
    \draw[color=mydark,->-=.5,line width=.5pt] plot [smooth,tension=.75]
    coordinates { (aa) (-2,-.25) (aaprime) };
    \node at (-3.6,.65) {\footnotesize $a$};
    \node at (-.35,.7) {\footnotesize $a'$};
    \node at (-2.25,-.35) {\footnotesize $p$};
    \node[circle,fill=mydark,inner sep=0pt,outer sep=0pt,minimum
    size=.75mm] (faa) at (3.25,.25) {};
    \node[circle,fill=mydark,inner sep=0pt,outer sep=0pt,minimum
    size=.75mm] (faaprime) at (4.75,-.25) {};
   \draw[color=mydark,->-=.5,line width=.5pt] (faa) to
   node[mylabel,auto,swap] {$f`p$} (faaprime);
   \draw[->,line width=.75,bend left] (.25,.75) to node[auto,mylabel]
   {$f$} (2.25,.75);
   \node at (3,.4) {\footnotesize $f a$};
   \node at (5,-.05) {\footnotesize $f a'$};
   \node at (-4.25,1.35) {$A$};
   \node at (5.25,1.25) {$B$};
  \end{tikzpicture}
  \caption{Representation of \texttt{maponpaths}.}
\end{figure}
We leave it as an exercise for the reader to verify that the
operation \verb|maponpaths| respects identity paths, as well as
composition and inverses of paths.

\section{Transport}\label{sec:transport}

Given a fibration $\pi:E\to B$ together with a path $f$ from $b$ to $b'$
in the base $B$, there is a continuous function $f_{!}:E_{b}\to
E_{b'}$ from the fiber $E_{b}$ of $\pi$ over $b$ to the fiber $E_{b'}$
over $b'$.  This operation $f_{!}$ of \emph{forward transport} along a
path is described in Coq as follows:
\begin{center}
  \begin{coqcode}
Definition transportf { B : UU } ( E : B -> UU ) { b b' : B } 
( f : paths b b' ) : E b -> E b'.
Proof.
  intros e. destruct f. assumption.
Defined.
\end{coqcode}
\end{center}
\begin{figure}[H]
  \centering
  \begin{tikzpicture}[color=mydark, fill=mylight, line width=1pt,scale=.75]
    \draw[fill=mylight] plot [smooth cycle,tension=.75] coordinates
    {(-2.5,0) (-2,1.5) (0,1.25) (1.75,1.5) (2.75,0) (1.75,-1) (-.5,0)};
    \draw[color=mydark,->-=.5,line width=.7pt] plot [smooth,tension=.5]
    coordinates { (-1.75,.5) (-1,.8)(1,0) (1.5,.25) };
    \node[circle,fill=mydark,inner sep=0pt,outer sep=4pt,minimum size=.75mm] (aa)
    at (-1.75,.5) {};
    \node[circle,fill=mydark,inner sep=0pt,outer sep=2pt,minimum size=.75mm] (bb)
    at (1.5,.25) {};
    \draw[dotted,line width=.6pt] (-1.75,3) to (aa);
    \draw[dotted,line width=.6pt] (1.5,3) to (bb);
    \draw[fill=mylight!50,line width=.8pt] plot [smooth cycle,tension=-.15] coordinates
    {(-2.5,2.5) (-2.5,4.5) (-1,4.5) (-1,2.5) };
    \draw[fill=mylight!50,line width=.8pt] plot [smooth cycle,tension=-.15] coordinates
    {(.75,2) (.75,4) (2.25,4) (2.25,2) };
    \node[circle,fill=mydark,inner sep=0pt,outer sep=2pt,minimum
    size=.75mm] (ee)
    at (-2,3.8) {};
    \node[circle,fill=mydark,inner sep=0pt,outer sep=2pt,minimum
    size=.75mm] (eef) at (1.5,3.5) {};
    \draw[->,dotted,line width=.6pt] (ee) to (eef);
    \node[font=\tiny] at (-2.15,3.9) {$\scriptstyle e$};
    \node[font=\tiny] at (1.75,3.7) {$\scriptstyle f_{!}(\!e)$};
    \node at (-1.95,.6) {$\scriptstyle b$};
    \node at (1.7,.4) {$\scriptstyle b'$};
    \node[font=\tiny] at (.25,.5) {$f$};
    \node[font=\tiny] at (-2.95,4.5) {$E(b)$};
    \node[font=\tiny] at (2.75,4) {$E(b')$};
    \node[font=\footnotesize] at (-2.55,1.45) {$B$};
  \end{tikzpicture} \label{fig:transportf}
  \caption{Forward transport.}
\end{figure}
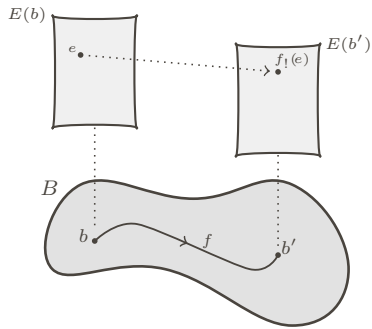
For a path $f$ as above, there is a corresponding operation
$f^{*}:E_{b'}\to E_{b}$ of \emph{backward transport} and it turns out
that $f_{!}$ and $f^{*}$ constitute a homotopy equivalence.

We will turn to briefly discuss homotopy and homotopy equivalence in
the setting of Coq before returning to forward and backward transport.

\subsection{Homotopy and homotopy equivalence}

Recall that for continuous functions $f,g:A\to B$, \emph{a homotopy from $f$
to $g$} is given by a continuous map $h:A\to B^{I}$ such that
\begin{center}
  \begin{tikzpicture}[auto]
    \node (UL) at (0,1.25) {$A$};
    \node (UR) at (2.5,1.25) {$B^{I}$};
    \node (B) at (1.25,0) {$B\times B$};
    \draw[->] (UL) to node[mylabel] {$h$} (UR);
    \draw[->,bend right=10pt] (UL) to node[mylabel,swap] {$\langle f,g\rangle$} (B);
    \draw[->,bend left=10pt] (UR) to (B);
  \end{tikzpicture}
\end{center}
commutes.  

In Coq, the type of homotopies between functions $f,g: A\to B$ is given by
\begin{center}
  \begin{coqcode}
Definition homot { A B : UU } ( f g : A -> B ) := forall x :A, paths ( f x ) ( g x ).
  \end{coqcode}
\end{center}
Here we encounter a new ingredient in Coq: the universal quantifier
\verb|forall|.  From the homotopical point of view, this
operation takes a fibration \verb|E : B -> UU| and gives back the
space \verb|forall x : B, E x| of all continuous sections of the
fibration.  That is, we should think of a point $s$ of this type as
corresponding to a continuous section
\begin{align*}
  \begin{tikzpicture}[auto]
    \node (UL) at (0,1.25) {$B$};
    \node (UR) at (2.5,1.25) {$E$};
    \node (B) at (1.25,0) {$B.$};
    \draw[->] (UL) to node[mylabel] {$s$} (UR);
    \draw[->,bend right=10pt] (UL) to node[mylabel,swap] {$1_{B}$} (B);
    \draw[->,bend left=10pt] (UR) to (B);
  \end{tikzpicture}
\end{align*}
One particular consequence of this is that if we are given a term 
\begin{center}
  \begin{coqcode}
s : ( forall x : B, E x )
  \end{coqcode}
\end{center}
and another term \verb|b : B|, then the term \verb|s| can be
\emph{applied} to the term \verb|b : B| to obtain a term of type 
\verb|E b|.  The result of applying \verb|s| to 
\verb|b| is denoted by
\begin{center}
  \begin{coqcode}
s b : E b.
  \end{coqcode}
\end{center}
We have more below to say about \verb|forall|.

Now, a map $f:A\to B$ is a \emph{homotopy equivalence} when there exists a
map $f':B\to A$ together with homotopies from $f'\circ f$ to $1_{A}$
and from $f\circ f'$ to $1_{B}$.  In this case, we say that $f'$ is a
\emph{homotopy inverse} of $f$.  Two spaces $A$ and $B$ are said to
have the same \emph{homotopy type} when there exists a homotopy
equivalence $f:A\to B$.

In Coq, we define the type of proofs that a map \verb|f : A -> B|
is a homotopy equivalence as follows:
\begin{center}
  \begin{coqcode}
Definition isheq { A B : UU } ( f : A -> B ) := total (fun f' : B -> A => dirprod (homot (funcomp f' f) (idfun _)) (homot (funcomp f f') (idfun _)) ).
  \end{coqcode}
\end{center}
Here it is worth pausing for a moment to consider the meaning of the
type \verb|isheq|.  Intuitively, \verb|isheq f| is the type
consisting of the data which one must provide in order to prove that
\verb|f| is a homotopy equivalence.  That is, a term of type
\verb|isheq f| consists of:
\begin{itemize}
\item a continuous map \verb|f' : B -> A|;
\item a homotopy from \verb|funcomp f' f| to the identity on \verb|B|;
\item a homotopy from \verb|funcomp f f'| to the identity on
  \verb|A|.
\end{itemize}
Indeed, by the definitions of \verb|total| and
\verb|dirprod| the terms of \verb|isheq f| can be regarded
as a tuple of such data.

\subsection{Forward and backward transport}

It turns out that, as mentioned above, the backward transport map
$f^{*}:E_{b'}\to E_{b}$ is a homotopy inverse of forward transport
$f_{!}$.  Denote by \verb|transportb| the backward transport
term.  It is often convenient to break the proofs of larger facts up
into smaller lemmas and we will do just this in order to show that
\verb|transportf E f| is a homotopy equivalence.  In particular,
we begin by proving that $f_{!}\circ f^{*}$ is homotopic to the
identity $1_{E_{b'}}$:
\begin{center}
  \begin{coqcode}
Lemma backandforth { B : UU } { E : B -> UU } { b b' : B } ( f : paths b b' ) ( e : E b' ) : homot ( funcomp ( transportb E f ) ( transportf E f ) ) ( idfun _ ).
Proof.
  intros x. destruct f. apply idpath. 
Defined.
  \end{coqcode}
\end{center}
Next, we prove that $f^{*}\circ f_{!}$ is homotopic to the identity
$1_{E_{b}}$ as \verb|forthandback| (we omit the proof because it
is identical to the proof of \verb|backandforth|):
\begin{center}
  \begin{coqcode}
Lemma forthandback { B : UU } { E : B -> UU } { b b' : B } ( f : paths b b' ) ( e : E b' ) : homot ( funcomp ( transportf E f ) ( transportb E f ) ) ( idfun _ ).
  \end{coqcode}
\end{center}
Using these lemmas we can finally prove that 
\verb|transportf E f| is a homotopy equivalence.
\begin{center}
  \begin{coqcode}
Lemma isheqtransportf { B : UU } ( E : B -> UU ) { b b' : B } ( f : paths b b' ) : isheq ( transportf E f ).
Proof.
  split with ( transportb E f ). split.
  apply backandforth. apply forthandback.
Defined.
  \end{coqcode}
\end{center}
\begin{figure}[H]
  \begin{tikzpicture}
    \node[smallcoqbox] (zero)  at (0,0) {%
      \begin{minipage}{5.2cm}
        \footnotesize
        \noindent\verb|1 subgoals, subgoal 1|

        ~

        \noindent\verb|B : UU|

        \noindent\verb|E : B -> UU|

        \noindent\verb|b : B|

        \noindent\verb|b' : B|
        
        \noindent\verb|f : paths b b'|

        \noindent\verb|============================|

        \noindent\verb| isheq (transportf E f)|

        \vphantom{\texttt{(transportf E f )) ( idfun ( E b' ) )}}
      \end{minipage}
    };
    \node[anchor=north east, inner sep=2pt] (titlezero) at
    (zero.north east) {\footnotesize\emph{Start of proof}};
    \node[smallcoqbox] (one) at (6,0) {%
      \begin{minipage}{5.2cm}
        \footnotesize
        \noindent\verb|2 subgoals, subgoal 1|
        
        ~
        
        \noindent\verb|B : UU|
        
        \noindent\verb|E : B -> UU|

        \noindent\verb|b : B|
        
        \noindent\verb|b' : B|
        
        \noindent\verb|f : paths b b'|
        
        \noindent\verb|============================|
        
        \noindent\verb| homot (funcomp (transportb E f)|
        
        \noindent\verb|  (transportf E f)) (idfun (E b'))|
      \end{minipage}
    };
    \node[anchor=north east, inner sep=2pt] (titleone) at
    (one.north east) {\footnotesize\emph{after} 
      \verb|split with; split|.};
    \node[smallcoqbox] (two)  at (0,-5.15) {%
      \begin{minipage}{5.2cm}
        \footnotesize
        \noindent\verb|1 subgoals, subgoal 1|

        ~

        \noindent\verb|B : UU|

        \noindent\verb|E : B -> UU|

        \noindent\verb|b : B|

        \noindent\verb|b' : B|
        
        \noindent\verb|f : paths b b'|

        \noindent\verb|============================|

        \noindent\verb| homot (funcomp (transportf E f)|
        
        \noindent\verb|  (transportb E f)) (idfun (E b))|
      \end{minipage}
    };
    \node[anchor=north east, inner sep=2pt] (titletwo) at
    (two.north east) {\footnotesize\emph{after} \verb|apply backandforth|};
  \end{tikzpicture}
  \caption{Coq output during the proof that forward transport is a
    homotopy equivalence.}
  \label{figure:isheqtransportf}
\end{figure}
There are several points to make about this proof.  The initial goal is
to supply a term of type \verb|isheq ( transportf E f )|.  Now,
this type is itself really of the form (you can see this in the proof
by entering \verb|unfold isheq|):
\begin{center}
  \begin{coqcode}
total (fun f' : E b' -> E b => dirprod (homot (funcomp f' (transportf
E f)) (idfun (E b'))) (homot (funcomp (transportf E f) f') (idfun (E b))))
  \end{coqcode}
\end{center}
and in general to construct a term of type \verb|total E|, for 
\verb|E : B -> UU|, it suffices (by virtue of the definition of 
\verb|total|) to give a term \verb|b| of type
\verb|B| together with a term of type \verb|E b|.  This is
captured in Coq by the command \verb|split with| and one should
think of \verb|split with b| as saying to Coq that you will
construct the required term using \verb|b| as the term of type
\verb|B| you are after.  Upon using this command, the goal will
automatically be updated to \verb|E b|.  In this case, 
entering \verb|split with ( transportb E f)| is the way to tell
Coq that we take \verb|transportb E f| to be the homotopy inverse
of \verb|transportf E f|.  So, after entering this
command the new goal becomes
\begin{center}
  \begin{coqcode}
dirprod 
 (homot (funcomp (transportb E f) (transportf E f)) (idfun (E b')))
 (homot (funcomp (transportf E f) (transportb E f)) (idfun (E b)))
  \end{coqcode}
\end{center}
As with \verb|total E|, in order to construct a term of type
\verb|dirprod A B| it suffices to supply terms of both types
\verb|A| and \verb|B|.  When given a goal of
the form \verb|dirprod A B|, we use the \verb|split| tactic
to tell Coq that we will supply separately the
terms of type \verb|A| and \verb|B| individually (as opposed
to providing a term by some other means).  (See Figure \ref{figure:isheqtransportf}
for the result of applying both \verb|split with| and
\verb|split| in the particular proof we are considering.)

The final new ingredient from the proof of \verb|isheqtransportf|
is the appearance of the tactic \verb|apply|.  When you have
proved a result in Coq and you are later
given a goal which is a (more or less direct) consequence of that
the result, then the tactic \verb|apply| will allow
you to apply the result.  In this case, the lemmas
\verb|backandforth| and \verb|forthandback| are exactly the
lemmas required in order to prove the remaining subgoals.

\subsection{Paths in the total space}

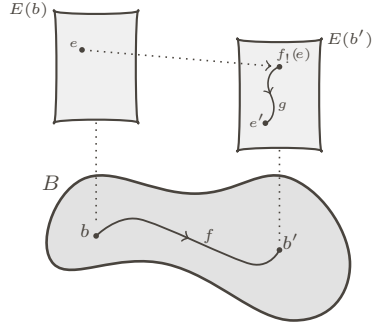
\begin{figure}[H]
  \centering
  \begin{tikzpicture}[color=mydark, fill=mylight, line width=1pt,scale=.75]
    \draw[fill=mylight] plot [smooth cycle,tension=.75] coordinates
    {(-2.5,0) (-2,1.5) (0,1.25) (1.75,1.5) (2.75,0) (1.75,-1) (-.5,0)};
    \draw[->-=.5,line width=.7pt] plot [smooth,tension=.5]
    coordinates { (-1.75,.5) (-1,.8)(1,0) (1.5,.25) };
    \node[circle,fill=mydark,inner sep=0pt,outer sep=4pt,minimum size=.75mm] (aa)
    at (-1.75,.5) {};
    \node[circle,fill=mydark,inner sep=0pt,outer sep=2pt,minimum size=.75mm] (bb)
    at (1.5,.25) {};
    \draw[dotted,line width=.6pt] (-1.75,3) to (aa);
    \draw[dotted,line width=.6pt] (1.5,3) to (bb);
    \draw[fill=mylight!50,line width=.8pt] plot [smooth cycle,tension=-.15] coordinates
    {(-2.5,2.5) (-2.5,4.5) (-1,4.5) (-1,2.5) };
    \draw[fill=mylight!50,line width=.8pt] plot [smooth cycle,tension=-.15] coordinates
    {(.75,2) (.75,4) (2.25,4) (2.25,2) };
    \draw[->-=.5,line width=.6pt] plot
      [smooth,tension=1] coordinates { (1.5,3.5) (1.3,3.25) (1.4,2.75) (1.25,2.5) };
    \node[circle,fill=mydark,inner sep=0pt,outer sep=2pt,minimum
    size=.75mm] (ee)
    at (-2,3.8) {};
    \node[circle,fill=mydark,inner sep=0pt,outer sep=2pt,minimum
    size=.75mm] (eef) at (1.5,3.5) {};
    \node[circle,fill=mydark,inner sep=0pt,outer sep=2pt,minimum
      size=.75mm] (eeprime) at (1.25,2.5) {};
    \draw[->,dotted,line width=.6pt] (ee) to (eef);
    \node[font=\tiny] at (-2.15,3.9) {$\scriptstyle e$};
    \node[font=\tiny] at (1.75,3.7) {$\scriptstyle f_{!}(\!e)$};
    \node[font=\tiny] at (1.55,2.8) {$\scriptstyle g$};
    \node[font=\tiny] at (1.1,2.55) {$\scriptstyle e'$};
    \node at (-1.95,.6) {$\scriptstyle b$};
    \node at (1.7,.4) {$\scriptstyle b'$};
    \node[font=\tiny] at (.25,.5) {$f$};
    \node[font=\tiny] at (-2.95,4.5) {$E(b)$};
    \node[font=\tiny] at (2.75,4) {$E(b')$};
    \node[font=\footnotesize] at (-2.55,1.45) {$B$};
  \end{tikzpicture}
  \caption{Paths in the total space.}
  \label{fig:pathintotalfiber}
\end{figure}

Using transport it is possible to give a complete characterization of
paths in the total space of a fibration \verb|E : B -> UU|.
Along these lines, the following lemma gives sufficient conditions for
the existence of a path in the total space:
\begin{center}
  \begin{coqcode}
Lemma pathintotalfiber { B : UU } { E : B -> UU } { x y : total E } ( f : paths ( pr1 x ) ( pr1 y ) ) ( g : paths ( transportf E f ( pr2 x ) ) ( pr2 y ) ) : paths x y.
Proof.
  intros. destruct x as [ x0 x1 ]. destruct y as [ y0 y1 ].  
  simpl in *. destruct f. destruct g. apply idpath.
Defined.
  \end{coqcode}
\end{center}
This lemma shows that, given points \verb|x| and \verb|y| of
the total space, in order to construct a path from \verb|x| to
\verb|y| it suffices to provide the following data:
\begin{itemize}
\item a path \verb|f| from \verb|pr1 x| to \verb|pr1 y|; and
\item a path \verb|g| from the result of transporting
  \verb|pr2 x| along \verb|f| to \verb|pr2 y|. 
\end{itemize}
This is illustrated in Figure \ref{fig:pathintotalfiber} in the
special case where \verb|x| is the pair \verb|pair b e| and
\verb|y| is the pair \verb|pair b' e'|.

Regarding the proof of \verb|pathintotalfiber|, it is worth
mentioning that here the effect of applying 
\verb|destruct x as [ x0 x1 ]| is that it tells Coq that we would like to consider the case where
\verb|x| is really of the form \verb|pair x0 x1|.  The only new
tactic here is \verb|simpl in *| which tells Coq to make any possible
simplifications to the terms appearing in the goal or hypotheses.
For example, in this case, Coq will simplify 
\verb|(pr1 (pair x0 x1))| to \verb|x0|.

On the other hand, if we are given a path \verb|f| from
\verb|x| to \verb|y| in the total space, there is an induced
path in the base given by 
\begin{center}
  \begin{coqcode}
Definition pathintotalfiberpr1 { B : UU } { E : B -> UU } { x y : total E } ( f : paths x y ) : paths ( pr1 x ) ( pr1 y ) := pr1 ` f.
  \end{coqcode}
\end{center}
Furthermore, we may transport \verb|pr2 x| along
\verb|pathintotalfiberpr1 f| and there is a path from the
resulting term to \verb|pr2 y|:
\begin{center}
  \begin{coqcode}
Definition pathintotalfiberpr2 { B : UU } { E : B -> UU } { x y : total E } ( f : paths x y ) : paths (transportf E ( pathintotalfiberpr1 f ) ( pr2 x )) ( pr2 y ).
Proof.
  intros. destruct f. apply idpath.
Defined.
  \end{coqcode}
\end{center}
Finally, we prove that every path in the total space is homotopic to
one obtained using \verb|pathintotalfiber|:
\begin{center}
  \begin{coqcode}
Lemma pathintotalfibercharacterization { B : UU } { E : B -> UU } { x
  y : total E } ( f : paths x y ) : paths f  (pathintotalfiber (pathintotalfiberpr1 f) (pathintotalfiberpr2 f) ).
Proof.
  intros. destruct f. destruct x as [ x0 x1 ]. apply idpath.
Defined.
  \end{coqcode}
\end{center}

\section{Weak equivalences and homotopy equivalences} \label{sec:bh}

In this section we introduce further basic homotopy theoretic notions
including the crucial notion of \emph{weak equivalence}.  Classically,
weak equivalences are those maps which induce isomorphisms on (all
higher) homotopy groups.  Between sufficiently nice spaces weak
equivalences turn out to coincide with homotopy equivalences and the
main goal of this section is to prove this fact in Coq.  The Coq
proofs in this section become increasingly sophisticated and
accordingly we will begin to pass quickly over the more basic features
of the proofs while at the same time giving their mathematical
interpretations.

\subsection{Contractibility}\label{sec:contr}

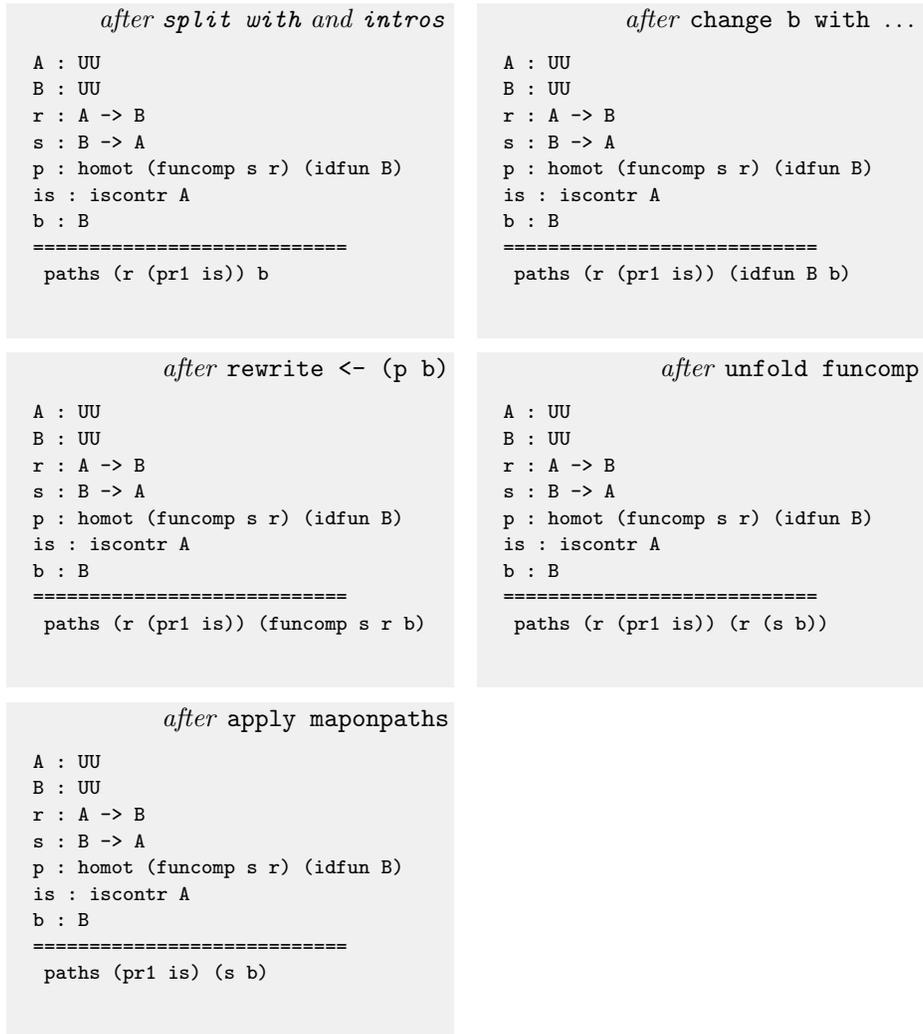
\begin{figure}[H]
  \begin{tikzpicture}
    \node[smallcoqbox] (zero)  at (0,0) {%
      \begin{minipage}{5.25cm}
        \footnotesize

        \noindent\verb|A : UU|

        \noindent\verb|B : UU|
        
        \noindent\verb|r : A -> B|

        \noindent\verb|s : B -> A|

        \noindent\verb|p : homot (funcomp s r) (idfun B)|
        
        \noindent\verb|is : iscontr A|

        \noindent\verb|b : B|

        \noindent\verb|============================|
        
        \noindent\verb| paths (r (pr1 is)) b|
      \end{minipage}
    };
    \node[anchor=north east, inner sep=2pt] (titlezero) at
    (zero.north east) {\emph{after \texttt{split with} and
        \texttt{intros}}};
    \node[smallcoqbox] (one) at (6.25,0) {%
      \begin{minipage}{5.25cm}
        \footnotesize

        \noindent\verb|A : UU|

        \noindent\verb|B : UU|
        
        \noindent\verb|r : A -> B|

        \noindent\verb|s : B -> A|

        \noindent\verb|p : homot (funcomp s r) (idfun B)|
        
        \noindent\verb|is : iscontr A|

        \noindent\verb|b : B|

        \noindent\verb|============================|
        
        \noindent\verb| paths (r (pr1 is)) (idfun B b)|
      \end{minipage}
    };
    \node[anchor=north east, inner sep=2pt] (titleone) at
    (one.north east) {\emph{after} \verb|change b with |\ldots};
    \node[smallcoqbox] (two) at (0,-4.65) {%
      \begin{minipage}{5.25cm}
        \footnotesize
        
        \noindent\verb|A : UU|

        \noindent\verb|B : UU|
        
        \noindent\verb|r : A -> B|

        \noindent\verb|s : B -> A|

        \noindent\verb|p : homot (funcomp s r) (idfun B)|
        
        \noindent\verb|is : iscontr A|

        \noindent\verb|b : B|

        \noindent\verb|============================|
        
        \noindent\verb| paths (r (pr1 is)) (funcomp s r b)|
      \end{minipage}
    };
    \node[anchor=north east, inner sep=2pt] (titletwo) at
    (two.north east) {\emph{after} \verb|rewrite <- (p b)|};
    \node[smallcoqbox] (three) at (6.25,-4.65) {%
      \begin{minipage}{5.25cm}
        \footnotesize
        
        \noindent\verb|A : UU|

        \noindent\verb|B : UU|
        
        \noindent\verb|r : A -> B|

        \noindent\verb|s : B -> A|

        \noindent\verb|p : homot (funcomp s r) (idfun B)|
        
        \noindent\verb|is : iscontr A|

        \noindent\verb|b : B|

        \noindent\verb|============================|
        
        \noindent\verb| paths (r (pr1 is)) (r (s b))|
      \end{minipage}
    };
    \node[anchor=north east, inner sep=2pt] (titlethree) at
    (three.north east) {\emph{after} \verb|unfold funcomp |};
    \node[smallcoqbox] (four) at (0,-9.3) {%
      \begin{minipage}{5.25cm}
        \footnotesize
        
        \noindent\verb|A : UU|

        \noindent\verb|B : UU|
        
        \noindent\verb|r : A -> B|

        \noindent\verb|s : B -> A|

        \noindent\verb|p : homot (funcomp s r) (idfun B)|
        
        \noindent\verb|is : iscontr A|

        \noindent\verb|b : B|

        \noindent\verb|============================|
        
        \noindent\verb| paths (pr1 is) (s b)|
      \end{minipage}
    };
    \node[anchor=north east, inner sep=2pt] (titlefour) at
    (four.north east) {\emph{after} \verb|apply maponpaths|};
  \end{tikzpicture}
  \label{fig:iscontrretract}
  \caption{Coq output during the proof of \texttt{iscontrretract}.}
\end{figure}

A space $A$ is said to be \emph{contractible} when the canonical map
$A\to 1$ to the one point space is a homotopy equivalence.  In Coq,
the notion of contractibility is captured by the following definition:
\begin{center}
  \begin{coqcode}
Definition iscontr ( A : UU ) := total ( fun center : A => forall b : A, paths center b ).
  \end{coqcode}
\end{center}
That is, a term of type \verb|iscontr A| consists of the data
required to prove that \verb|A| is contractible:
\begin{itemize}
\item a point \verb|center| of \verb|A|, which we will
  sometimes refer to as the \emph{center of contraction}; and
\item a continuous assignment of, to each \verb|b : A|, a path
  from \verb|center| to \verb|b|.
\end{itemize}
There are various facts about contractible space (e.g., $1$ is
contractible, contractible spaces satisfy the principle of ``proof
irrelevance'', and so forth) which are ready consequences of this
definition.  However, we will leave the investigation of such matters
to the reader and merely include one important fact about contractible
spaces: homotopy retracts of contractible spaces are contractible.  That is,
given continuous functions $r : A \to B$ and $s: B\to A$ together with
a homotopy $p$ as indicated in the following diagram:
\begin{align}\label{eq:hretract}
  \begin{minipage}{.25\linewidth}
    \begin{tikzpicture}[auto]
    \node (UL) at (0,1.5) {$B$};
    \node (UR) at (3,1.5) {$A$};
    \node (LC) at (1.5,0) {$B$};
    \draw[->] (UL) to node[mylabel] {$s$} (UR);
    \draw[->,bend left] (UR) to node[mylabel] {$r$} (LC);
    \draw[->,bend right] (UL) to node[mylabel,swap] {$1_{B}$} (LC);
    \node (ppp) at (1.5,.75) {$\Downarrow$}; 
    \node (pp) at (1.75,.75) {\footnotesize $p$};
  \end{tikzpicture}
  \end{minipage}
\end{align}
if the space $A$ is contractible, then so is the space $B$.  Note that
in the diagram (\ref{eq:hretract}) we employ the convention, familiar
from higher-dimensional category theory, of indicating the homotopy $p$
by a double arrow $\Rightarrow$.  In Coq:
\begin{center}
  \begin{coqcode}
Lemma iscontrretract { A B : UU } { r : A -> B } { s : B -> A } (p : homot ( funcomp s r ) (idfun _)) (is : iscontr A) : iscontr B.
Proof.
  split with ( r ( pr1 is ) ).  intros b. 
  change b with ( idfun B b ). rewrite <- ( p b ). unfold funcomp. 
  apply maponpaths. apply ( pr2 is ).
Defined.
  \end{coqcode}
\end{center}
The output of Coq during the proof can be found in Figure
\ref{fig:iscontrretract}.  First, we need to provide a center of contraction for \verb|B|.
The center of contraction for \verb|A| is the term 
\verb|pr1 is|. By entering 
\verb|split with ( r (pr1 is ) )| we tell Coq to take
the center of contraction of \verb|B| to be the term
\verb|r ( pr1 is )|.  After this, the goal becomes
\begin{center}
  \begin{coqcode}
forall b : B, paths (r (pr1 is)) b
  \end{coqcode}
\end{center}
which is to say that we must prove that there is a path from
\verb|r (pr1 is)| to any other point of \verb|B|.  Next,
after using \verb|intros b|, we enter the command 
\begin{center}
  \begin{coqcode}
change b with ( idfun B b )
  \end{coqcode}
\end{center}
to replace the term \verb|b| in the goal with the \emph{equal}
term
\verb|idfun B b|.  Observe that the term \verb|p b| has type
\begin{center}
  \begin{coqcode}
paths (funcomp s r b) (idfun B b).
  \end{coqcode}
\end{center}
The tactic \verb|rewrite| allows us to replace each occurrence of
a term in the goal by a term which is in the same path component.
In this case, because \verb|p b| is a path from
\verb|funcomp s r b| to \verb|idfun B b|, the result of
rewriting with \verb|rewrite <- (p b)| is to replace each
occurrence of \verb|idfun B b| in the goal by 
\verb|funcomp s r b|.  Here the backward arrow \verb|<-|
indicates that we are rewriting the path \verb|p b| from right to
left.

Next, we use \verb|unfold funcomp| to replace each occurrence of 
\verb|funcomp| in the goal by its definition.  Then
\verb|apply maponpaths| tells Coq that we will use 
\verb|maponpaths| to construct the goal.  Because Coq recognizes
that the hypotheses of \verb|maponpaths| require us to provide a
path the effect of this is to remove both applications of the function
\verb|r| on both the right and left.  Finally, the required path
can be constructed using \verb|pr2 is|.

It is an immediate consequence of \verb|iscontrretract| that
spaces which are homotopy equivalent to contractible spaces are also
contractible.  This fact is captured by the definitions
\begin{center}
  \begin{coqcode}
Definition iscontrandheq { A B : UU } { f : A -> B } ( p : isheq f ) (is : iscontr A) : iscontr B := iscontrretract (pr1 (pr2 p)) is.
  \end{coqcode}
\end{center}
and
\begin{center}
  \begin{coqcode}
Definition iscontrandheqinv { A B : UU } { f : A -> B } ( p : isheq f ) ( is : iscontr B ) : iscontr A := iscontrretract ( pr2 ( pr2 p ) ) is. 
  \end{coqcode}
\end{center}

\subsection{Homotopy fibers}

Given a map \verb|f : A -> B| and a point \verb|b| of 
\verb|B|, we define the \emph{homotopy fiber of \texttt{f}
  over \texttt{b}} as follows:
\begin{center}
  \begin{coqcode}
Definition hfiber { A B : UU } ( f : A -> B ) ( b : B ) : UU := 
  total ( fun x : A => paths ( f x ) b ).
  \end{coqcode}
\end{center}
The homotopy fiber of $f$ over $b$ is the homotopical analogue of the
ordinary fiber $f^{-1}(b)$.  So, a typical point of the homotopy fiber
is a pair consisting of a point $a$ of $A$ together with a path 
\begin{align*}
  \begin{tikzpicture}
    \node (A) at (0,0) {$f(a)$};
    \node (B) at (2,0) {$b$};
    \draw[->-=.5] (A) to (B);
  \end{tikzpicture}
\end{align*}
in $B$.

There are various reasons for considering homotopy fibers.  Homotopy
fibers play a role for fibrations analogous to that played by ordinary
fibers for fiber bundles.  For us, the interest in homotopy fibers
comes from their presence in the definition of weak equivalences given below.

From the point of view of category theory, the homotopy fiber of $f$
over $b$ is the $\infty$-groupoid version of the comma category
$(f\downarrow b)$.  As in the categorical setting, there are actions
of the operation of taking the homotopy fiber over a fixed point on
maps over the base $B$.  That is, given
continuous maps $f:A\to B$, $g:C\to B$ and $h:A\to C$ together with a
homotopy $p$ from $g\circ h$ to $f$, there is an induced map 
\begin{align*}
  \texttt{hfiber}\;f\;b \to \texttt{hfiber}\;g\;b.
\end{align*}
In Coq, this map is defined as 
\begin{center}
  \begin{coqcode}
Definition hfiberact { A B C : UU } { f : A -> B } { g : C -> B } ( h : A -> C ) ( p : homot ( funcomp h g ) f ) ( b : B ) : hfiber f b -> hfiber g b 
:= fun a => pair (h ( pr1 a )) (pathscomp (p ( pr1 a )) ( pr2 a )).
  \end{coqcode}
\end{center}
That is, \verb|hfiberact h p b| sends an element
\begin{align*}
  \begin{tikzpicture}[auto]
    \node (A) at (0,0) {$f(a)$};
    \node (B) at (2,0) {$b$};
    \draw[->-=.5] (A) to node[mylabel] {$i$} (B);
  \end{tikzpicture}
\end{align*}
of $\texttt{hfiber}\;f\;b$ to the point of $\texttt{hfiber}\;g\;b$
given by the composite path
\begin{align*}
  \begin{tikzpicture}[auto]
    \node (A) at (0,0) {$g\bigl(h(a)\bigr)$};
    \node (B) at (2,0) {$f(a)$};
    \node (C) at (4,0) {$b$.};
    \draw[->-=.5] (A) to node[mylabel] {$p(a)$} (B);
    \draw[->-=.5] (B) to node[mylabel] {$i$} (C);
  \end{tikzpicture}
\end{align*}
As a special case of this construction, when we are given a continuous
maps $f:A\to B$ and $g:B\to C$ together with a point $c$ of $C$, there
is an induced map 
\begin{align*}
\texttt{hfiber}\;(g\circ f)\;c\to\texttt{hfiber}\;g\;c
\end{align*}
obtained by applying 
\verb|hfiberact| with the identity homotopy $1_{g\circ f}$.
In Coq, this is obtained as the following definition:
\begin{center}
  \begin{coqcode}
Definition maponhfiber { A B C : UU } ( f : A -> B ) ( g : B -> C ) ( c : C ) : hfiber ( funcomp f g ) c -> hfiber g c := fun a => pair ( f ( pr1 a ) ) ( pr2 a ).
  \end{coqcode}
\end{center}
In the case where we have a homotopy retract as in (\ref{eq:hretract})
above, we obtain, a further map as follows:
\begin{center}
  \begin{coqcode}
Definition secmaponhfiber { A B : UU } {r : A -> B} { s : B -> A } ( p : homot ( funcomp s r ) ( idfun _ ) ) ( a : A ) : hfiber s a -> hfiber ( funcomp r s ) a := fun b => pair (s ( pr1 b)) (pathscomp (s ` ( p ( pr1 b ) )) ( pr2 b )).
  \end{coqcode}
\end{center}
That is, \verb|secmaponhfiber| sends
\begin{align*}
  \begin{tikzpicture}[auto]
    \node (A) at (0,0) {$s(b)$};
    \node (B) at (2,0) {$a$};
    \draw[->-=.5] (A) to node[mylabel] {$i$} (B);
  \end{tikzpicture}
\end{align*}
in $\texttt{hfiber}\;s\;a$ to the term in $\texttt{hfiber}\;(s\circ
r)\;a$ given by the composite path
\begin{align*}
  \begin{tikzpicture}[auto]
    \node (A) at (0,0) {$s\bigl(r(s(b))\bigr)$};
    \node (B) at (3,0) {$s(b)$};
    \node (C) at (5,0) {$a$.};
    \draw[->-=.5] (A) to node[mylabel] {$s`p(b)$} (B);
    \draw[->-=.5] (B) to node[mylabel] {$i$} (C);
  \end{tikzpicture}
\end{align*}
It turns out that, in this situation, $\texttt{hfiber}\;s\;a$ is a
homotopy retract of $\texttt{hfiber}\;(s\circ r)\;a$:
\begin{align*}
    \begin{tikzpicture}[auto]
    \node (UL) at (0,1.5) {$\texttt{hfiber}\;s\;a$};
    \node (UR) at (4.5,1.5) {$\texttt{hfiber}\;(s\circ r)\;a$};
    \node (LC) at (1.5,0) {$\texttt{hfiber}\;s\;a$};
    \draw[->] (UL) to node[mylabel] {$\texttt{secmaponhfiber}\;p\;a$} (UR);
    \draw[->,bend left] (UR) to node[mylabel] {$\texttt{maponhfiber}\;r\;s\;a$} (LC);
    \draw[->,bend right] (UL) to node[mylabel,swap] {$1_{\texttt{hfiber}\;s\;a}$} (LC);
    \node (ppp) at (1.5,.75) {$\Downarrow$};
  \end{tikzpicture}
\end{align*}
Here, as in (\ref{eq:hretract}), we indicate the existence of a
homotopy in this diagram by a double arrow $\Rightarrow$.  In Coq,
this fact is captured by the following lemma:
\begin{center}
  \begin{coqcode}
Lemma secmaponhfiberissec {A B : UU} { r : A -> B } { s : B -> A } ( p : homot ( funcomp s r ) ( idfun _ ) ) ( a : A ) : homot (funcomp  ( secmaponhfiber p a ) ( maponhfiber r s a )) (idfun _).
Proof.
 intros b. destruct b as [ b i ]. unfold funcomp, idfun in *. simpl. 
 apply (@pathintotalfiber _ _ (pair (r (s b)) _) (pair b i) (p b)).
 rewrite transportfandhfiber. unfold secmaponhfiber. simpl.
 rewrite <- isassocpathscomp. rewrite islinvpathsinv. apply idpath. 
Defined.
  \end{coqcode}
\end{center}
First, after the first application of the \verb|simpl| tactic we
have simplified the goal and hypotheses to the extent that we are now
in the following situation:
\begin{center}
  \begin{coqcode}
  A : UU
  B : UU
  r : A -> B
  s : B -> A
  p : homot (fun x : B => r (s x)) (fun x : B => x)
  a : A
  b : B
  i : paths (s b) a
  ============================
   paths (maponhfiber r s a (secmaponhfiber p a (pair b i))) 
     (pair b i)
  \end{coqcode}
\end{center}
Because we are asked to construct a path in a total space we will make
use of \verb|pathintotalfiber|.  However, in order to help Coq
understand exactly what we are trying to do we must provide some of
the implicit arguments of \verb|pathintotalfiber| explicitly.
This is accomplished here using the symbol \verb|@|.  The
underscores are the arguments which we leave for Coq to guess and the
other arguments are those we are explicitly providing for Coq.  That
is, we are telling Coq that in order to achieve the current goal it
suffices to construct using \verb|pathintotalfiber| a path from
\verb|pair ( r ( s b ) ) _ | to \verb|pair b i| where the
path from \verb|r ( s b )| to \verb|b| we use is 
\verb|p b|.  After this, the goal becomes:
\begin{center}
  \begin{coqcode}
paths
 (transportf (fun x : B => paths (s x) a) (p b)
   (pr2 (pair (r (s b)) (pr2 (secmaponhfiber p a (pair b i))))))
 (pr2 (pair b i))
  \end{coqcode}
\end{center}
Here, as is often the case, it is convenient to know that the result of applying forward
transport can be decomposed in a specific way.  In this case, it turns
out that, for \verb|j| a path from \verb|v| to \verb|w|
in \verb|B| and a path \verb|k : paths ( s v ) a|,
\begin{center}
  \begin{coqcode}
transportf ( fun x : B => paths ( s x ) a ) j k
  \end{coqcode}
\end{center}
is homotopic (in the sense of path homotopy) to the following
composite path:
\begin{align*}
  \begin{tikzpicture}[auto]
    \node (UL) at (0,0) {$s w$};
    \node (UR) at (2,0) {$s v$};
    \node (RR) at (4,0) {$a$.};
    \draw[->-=.5] (UL) to node[mylabel] {$(s ` j )^{-1}$} (UR);
    \draw[->-=.5] (UR) to node[mylabel] {$k$} (RR);
  \end{tikzpicture}
\end{align*}
This fact is captured by the lemma \verb|transportfandhfiber|
which has a trivial proof that we omit.  Returning to the proof of \verb|secmaponhfiberissec|, we rewrite
using \verb|transportfandhfiber| and then simplify using 
\verb|unfold| and \verb|simpl| to arrive at the goal
\begin{center}
  \begin{coqcode}
paths (pathscomp (pathsinv (s ` p b)) (pathscomp (s ` p b) i)) i
  \end{coqcode}
\end{center}
which is an immediate consequence of associativity of path composition
and the fact that \verb|pathsinv| is an inverse with respect to
path composition.

In addition to acting on maps via \verb|hfiberact|, there is also
an action of the homotopy fiber operation on homotopies.  Namely, a
homotopy from a map $f:A\to B$ to a map $g:A\to B$ induces, for $b :
B$, a map $\texttt{hfiber}\;f\;b\to\texttt{hfiber}\;g\;b$ as follows:
\begin{center}
\begin{coqcode}
Definition hfiberandhomot { A B : UU } { f g : A -> B } ( b : B ) ( p : homot f g ) : hfiber f b -> hfiber g b := fun a => pair (pr1 a) (pathscomp (pathsinv (p (pr1 a ))) (pr2 a)). 
\end{coqcode}
\end{center}
Similarly, there is a corresponding map
\begin{align*}
  \texttt{hfiber}\;g\;b\to\texttt{hfiber}\;f\;b
\end{align*}
because the homotopy relation is symmetric.
\begin{center}
\begin{coqcode}
Definition hfiberandhomotinv { A B : UU } { f g : A -> B } ( b : B ) ( p : homot f g ) : hfiber g b -> hfiber f b := fun a => pair ( pr1 a ) (pathscomp ( ( p ( pr1 a ) ) ) ( pr2 a )). 
\end{coqcode}
\end{center}
Finally, these two maps constitute a homotopy equivalence of spaces as
the following lemma confirms:
\begin{center}
\begin{coqcode}
Lemma isheqhfiberandhomot { A B : UU } { f g : A -> B } ( b : B ) ( p : homot f g ) : isheq ( hfiberandhomot b p ).
\end{coqcode}
\end{center}
We leave the proof of \verb|isheqhfiberandhomot| as an exercise
for the reader (the proof can also be found in the accompanying Coq
file). 

\subsection{Weak equivalences}\label{section:weq_def}

Classically, a map $f:A\to B$ is a \emph{weak equivalence} when it
induces isomorphisms on homotopy groups.  However, we will give a
definition, following Voevodsky \cite{Vo2012a}, which makes sense without referring
to homotopy groups.  In Section \ref{section:gradth} below we will
show that these weak equivalences coincide, for the ``nice spaces'' we
are considering, with the homotopy equivalences (and therefore also with
the classical weak equivalences).  We start with the definition:
\begin{center}
  \begin{coqcode}
Definition isweq { A B : UU } ( f : A -> B ) := ( forall b : B, iscontr ( hfiber f b ) ).
  \end{coqcode}
\end{center}
From the ``propositions as types'' point of view, the weak
equivalences correspond to the bijections.  The space of weak
equivalences from $A$ to $B$ is defined as follows:
\begin{center}
  \begin{coqcode}
Definition weq ( A B : UU ) := total ( fun f : A -> B => isweq f ).
  \end{coqcode}
\end{center}
That is, a typical term of type \verb|weq A B| is a pair
consisting of a map \verb|f : A -> B| together with a term of
type \verb|isweq f|.  The most basic example of a weak equivalence
is the identity function \verb|idfun A : A -> A|.  It is
straightforward to construct the required proof that this is a weak
equivalence.  We denote by \verb|isweqidfun A| this term of type 
\verb|isweq ( idfun A )| and we then adopt the following definition:
\begin{center}
  \begin{coqcode}
Definition idweq ( A : UU ) := pair ( idfun A ) ( isweqidfun A ).
  \end{coqcode}
\end{center}
That is, \verb|idweq A| is the representative of the identity
function on \verb|A| as a weak equivalence.  

Given a map \verb|f : A -> B| together with a proof 
\verb|is : isweq f| that it is a weak equivalence, if 
we are given a point \verb|b : B|, then there is a corresponding 
center of contraction in the homotopy fiber of \verb|f| over 
\verb|b| given explicitly by the term \verb|pr1 ( is b )|.
Because it is often convenient to make use the actual term of type 
\verb|A| corresponding to this center of contraction we introduce
the following nomenclature:
\begin{center}
  \begin{coqcode}
Definition weqpreimage { A B : UU } { f : A -> B } ( is : isweq f ) ( b : B ) : A := pr1 ( pr1 ( is b ) ).
  \end{coqcode}
\end{center}
That is, \verb|weqpreimage is b| is the (homotopically) canonical
element in the homotopy fiber of \verb|f| over \verb|b|.
Similarly, we name the path from \verb|f ( weqpreimage is b )| to 
\verb|b| as follows:
\begin{center}
  \begin{coqcode}
Definition weqpreimageeq { A B : UU } { f : A -> B } ( is : isweq f ) ( b : B ) : paths ( f ( weqpreimage is b ) ) b := pr2 ( pr1 ( is b ) ).
  \end{coqcode}
\end{center}
That is, we have:
\begin{align*}
   \begin{tikzpicture}[auto]
    \node (A) at (0,0) {$f(\texttt{weqpreimage is }b)$};
    \node (B) at (5,0) {$b$.};
    \draw[->-=.5] (A) to node[mylabel] {\texttt{weqpreimageeq is }$b$} (B);
  \end{tikzpicture}
\end{align*}
It is also clear from the definition of \verb|weqpreimage| that
if we are given any point \verb|a : A| and path \verb|g|
from \verb|f a| to \verb|b|, there exists a path
\begin{align*}
 \begin{tikzpicture}[auto]
    \node (A) at (0,0) {$\texttt{weqpreimage is }b$};
    \node (B) at (5,0) {$a$.};
    \draw[->-=.5] (A) to node[mylabel] {\texttt{weqpreimageump1 is
      }$b$ $g$} (B);
  \end{tikzpicture}
\end{align*}
We leave the construction of the term \verb|weqpreimageump1| to
the reader.  Finally, the operation of taking the preimage is
injective in the sense that if there is a path
\begin{align*}
  \begin{tikzpicture}[auto]
    \node (A) at (0,0) {$\texttt{weqpreimage is }b$};
    \node (B) at (5,0) {$\texttt{weqpreimage is }b'$,};
    \draw[->-=.5] (A) to node[mylabel] {$g$} (B);
  \end{tikzpicture}
\end{align*}
then there exists an induced path \verb|weqpreimageump2 is g|
from \verb|b| to \verb|b'|.

\subsection{Weak equivalences and homotopy equivalences}\label{section:gradth}

Given a weak equivalence \verb|f : A -> B|, we construct a
homotopy inverse of \verb|f| as follows:
\begin{center}
  \begin{coqcode}
Definition weqinv { A B : UU } ( f : weq A B ) : B -> A := fun x => weqpreimage ( pr2 f ) x.
  \end{coqcode}
\end{center}
Using the observations from Section \ref{section:weq_def} it is
straightforward to prove the following lemmas:
\begin{center}
  \begin{coqcode}
Lemma weqinvislinv { A B : UU } ( f : weq A B ) : homot ( funcomp ( weqinv f ) ( pr1 f ) ) ( idfun _ ).
  \end{coqcode}
\end{center}
and
\begin{center}
  \begin{coqcode}
Lemma weqinvisrinv { A B : UU } ( f : weq A B ) : homot ( funcomp ( pr1 f ) ( weqinv f ) ) ( idfun _ ).
  \end{coqcode}
\end{center}
which together constitute the proof of:
\begin{center}
  \begin{coqcode}
Lemma weqtoheq { A B : UU } { f : A -> B } (is : isweq f) : isheq f.
  \end{coqcode}
\end{center}
Classically it should be noted that the analogue of 
\verb|weqtoheq| only applies in the case of reasonably nice
spaces.  Interestingly, the proof of the converse of \verb|weqtoheq| is not
entirely trivial.  (Categorically, the proof of the converse makes
crucial use of the fact that categorical equivalences can be
transformed into adjoint equivalences \cite{MacLane:1971tv}.)  

The proof of the converse of \verb|weqtoheq| can nicely be broken
up into two lemmas.  The first lemma states that, for a homotopy
retract as in (\ref{eq:hretract}) above, if the homotopy fiber of the
composite $s\circ r$ over a point $a$ of $A$ is contractible, then the
homotopy fiber of $s$ over $a$ is also contractible:
\begin{center}
  \begin{coqcode}
Lemma iscontrhfiberandhretract { A B : UU } { r : A -> B } { s : B -> A } (p : homot ( funcomp s r ) ( idfun _ )) ( a : A ) : iscontr ( hfiber ( funcomp r s ) a ) -> iscontr ( hfiber s a ).
Proof.
 intros is. 
 apply ( @iscontrretract ( hfiber ( funcomp r s ) a ) ( hfiber s a ) 
  ( maponhfiber _ _ _ ) ( secmaponhfiber p a ) ). 
 apply secmaponhfiberissec. assumption. 
Defined.    
  \end{coqcode}
\end{center}
As the Coq code shows, the proof of this fact is an immediate
consequence of \verb|iscontrretract| and
\verb|secmaponhfiberissec|.

The second lemma states that if there is a homotopy from $f:A\to B$ to
$f':A\to B$
and the homotopy fiber of $f'$ over a point $b : B$ is contractible,
then the homotopy fiber of $f$ over $b$ is also contractible:
\begin{center}
  \begin{coqcode}
Lemma iscontrhfiberandhomot { A B : UU } { f f' : A -> B } ( h : homot
f f' ) ( b : B ) : iscontr ( hfiber f' b ) -> iscontr ( hfiber f b ).
Proof.
  intros is. apply ( iscontrandheqinv ( isheqhfiberandhomot b h ) ). 
  assumption.
Defined.
  \end{coqcode}
\end{center}
In this case the proof is an immediate consequence of the fact 
(\verb|iscontrandheqinv|) that being contractible is preserved by
homotopy equivalences together with the fact
(\verb|isheqhfiberandhomot|) that homotopies induce homotopy
equivalences of homotopy fibers.

Using these two lemmas we obtain the converse of \verb|weqtoheq|,
called following Voevodsky the \emph{Grad Theorem}, as follows:
\begin{center}
\begin{coqcode}
Theorem gradth { A B : UU } { f : A -> B } ( is : isheq f ) : isweq f.
Proof.
  intro b. destruct is as [ f' is ]. 
  apply ( @iscontrhfiberandhretract B A f' f ( pr2 is ) ). 
  apply (@iscontrhfiberandhomot B _ (funcomp f' f) (idfun B) (pr1 is)). 
  apply isweqidfun.
Defined.
\end{coqcode}
\end{center}
After applying \verb|intro b| and \verb|destruct| we find
ourselves in the following situation:
\begin{center}
  \begin{coqcode}
  A : UU
  B : UU
  f : A -> B
  f' : B -> A
  is : dirprod (homot (funcomp f' f) (idfun B))
         (homot (funcomp f f') (idfun A))
  b : B
  ============================
   iscontr (hfiber f b)
  \end{coqcode}
\end{center}
Because \verb|f'| is the homotopy inverse of \verb|f| it
then suffices, by \verb|iscontrhfiberandhretract| to show that
the composite $f\circ f'$ has a contractible homotopy fiber over $b$.
This is the effect of the first \verb|apply|, after which the Coq
output is
\begin{center}
  \begin{coqcode}
  A : UU
  B : UU
  f : A -> B
  f' : B -> A
  is : dirprod (homot (funcomp f' f) (idfun B))
         (homot (funcomp f f') (idfun A))
  b : B
  ============================
   iscontr (hfiber (funcomp f' f) b)
  \end{coqcode}
\end{center}
We now observe that there is a homotopy from $f\circ f'$ to the
identity function on $B$ and therefore after applying
\verb|iscontrhfiberandhomot| it remains only to prove that the
identity on $B$ is a weak equivalences, which is precisely the content
of \verb|isweqidfun|.

\section{The Univalence Axiom and some consequences}\label{sec:hlevels}

In this section we will describe a number of constructions and results
which are more closely related to the univalent approach.  Because it
would take us to far afield in an introductory paper such as this, we
will merely mention a number of the results and display some of the
corresponding Coq code.  That is to say, the development given here is
not self contained: there are many definitions and lemmas we do not
give that would be required in order to obtain all of the results
described here.

\subsection{An alternative characterization of the Univalence Axiom}

Before giving the explicit statement of the Univalence Axiom, we will
first require a map which turns a path in the universe \verb|UU|
into a weak equivalence.  For a change of pace, we give a direct
definition of this map as follows:
\begin{center}
  \begin{coqcode}
Definition eqweqmap { A B : UU } ( p : paths A B ) : weq A B 
  := match p with idpath => ( idweq _ ) end.
  \end{coqcode}
\end{center}
That is, \verb|eqweqmap| is the map from the path space
\verb|paths A B| to the space \verb|weq A B| of weak
equivalences induced by the operation of sending the identity path on
\verb|A| to the identity weak equivalence \verb|idweq A| on 
\verb|A|.

We then define the type 
\begin{center}
  \begin{coqcode}
Definition isweqeqweqmap := forall A B : UU, isweq (@eqweqmap A B).
  \end{coqcode}
\end{center}
The \emph{Univalence Axiom} then states that there is a term of type
\verb|isweqeqweqmap|.  However, before explicitly adding the
Univalence Axiom as an assumption, we would first like to give a
logically equivalent version of this principle.  The idea behind
this equivalent form of the Univalence Axiom, which can be seen easily
by considering the semantic version of the Univalence Axiom (i.e.,
what the axiom says in terms of models), states that the type of weak
equivalences is inductively generated by the terms of the form
\verb|idweq A|.
In particular, they are inductively generated by these terms in the
same way that the path space construction \verb|paths| is
inductively generated by the identity paths.  Formally, we define the
type
\begin{center}
  \begin{coqcode}
Definition weqindelim := 
forall E : 
total ( fun x : UU => total ( fun y : UU => weq x y ) ) -> UU,
forall p : ( forall x : UU, E ( pair x ( pair x ( idweq x ) ) ) ), 
forall x y : UU, forall z : weq x y, E ( pair x ( pair y z ) ).
  \end{coqcode}
\end{center}
Intuitively, let $B$ be the space with points given by the following
data:
\begin{itemize}
\item small spaces $X$ and $Y$;
\item weak equivalences $f$ from $X$ to $Y$.
\end{itemize}
Then, for a fibration $E\to B$ over $B$, if there exists a proof $p$ (term)
that each fiber $E_{1_{X}}$ over identity weak equivalences $1_{X}$ is
inhabited, then a term of type $\texttt{weqindelim}\;E\;p$ yields a proof
that \emph{every} fiber $E_{f}$ is inhabited.  Thinking of $E$ as a
``property'' of weak equivalences, this states that in order to prove
that a property (definable type theoretically) of weak equivalences holds it
suffices to prove that the property holds of identity weak
equivalences.  We refer to this as \emph{induction on weak equivalences}.

In order to prove that if the Univalence Axiom holds, then the
induction principle \verb|weqindelim| also holds we make use of
the following lemma, the proof of which is immediate.
\begin{center}
  \begin{coqcode}
Lemma weqind0 
( E : total (fun x : UU => total ( fun y : UU => weq x y ) ) -> UU) 
( p : ( forall x : UU, E ( pair x ( pair x ( idweq x ) ) ) ) ) : 
(forall x y : UU,
  forall z : paths x y, E ( pair x ( pair y ( eqweqmap z ) ) )).    
  \end{coqcode}
\end{center}
This lemma states that the induction principle in question holds for
weak equivalences of the form $\texttt{eqweqmap}(f)$ for $f$ a path
between small spaces.  Using this, we obtain the following:
\begin{center}
  \begin{coqcode}
Definition weqind ( univ : isweqeqweqmap ) : weqindelim := 
fun E p A B f => 
 transportf ( fun z => E ( pair A ( pair B z ) ) )
  ( weqeqmaplinv univ f ) ( weqind0 E p A B (weqeqmap univ f) ). 
  \end{coqcode}
\end{center}
Here we are assuming that the Univalence Axiom holds.
I.e., that there is a term \verb|univ| of type
\verb|isweqeqweqmap|.  Then, given all of the data of the
induction principle, in order to obtain a proof that the fiber $E_{f}$
is inhabited we observe that, by \verb|weqind0|, there exists a
term $e$ in the 
\begin{align*}
  E_{\texttt{eqweqmap}(\texttt{weqeqmap}(f))},
\end{align*}
where \verb|weqeqmap| is the homotopy inverse of 
\verb|eqweqmap|, which exists by the fact that we are assuming
the Univalence Axiom.  It then suffices to transfer $e$ along the path
from $\texttt{eqweqmap}(\texttt{weqeqmap}(f))$ to $f$ (here called
\verb|weqeqmaplinv|). 

We note that, when the Univalence Axiom holds, the following
computation principle corresponding to \verb|weqindelim| also holds:
\begin{center}
  \begin{coqcode}
Definition weqindcomp ( rec : weqindelim ) := 
forall  E :
 total ( fun x : UU => total ( fun y : UU => weq x y ) ) -> UU, 
forall p : ( forall x : UU, E ( pair x ( pair x ( idweq x ) ) ) ), 
forall x : UU, paths ( rec E p x x ( idweq x ) ) ( p x ).
  \end{coqcode}
\end{center}
Explicitly, we have the following Theorem:
\begin{center}
  \begin{coqcode}
Theorem weqcomp (univ : isweqeqweqmap) : weqindcomp ( weqind univ ).
  \end{coqcode}
\end{center}
The proof of this is slightly involved and we leave the details to the
reader (they can also be found in the companion Coq file for this
tutorial).  We refer to this computation principle as the
\emph{computation principle for weak equivalences}.

It turns out that the converse implications also hold: in order to
prove that the Univalence Axiom holds, it suffices to show that there
are terms \verb|rec : weqindelim| and 
\verb|reccomp : weqindcomp rec|.  That is, we have the following
Theorem:
\begin{center}
\begin{coqcode}
Theorem univfromind {rec : weqindelim} ( reccomp : weqindcomp rec ) 
 : isweqeqweqmap. 
\end{coqcode}
\end{center}
In order to prove this, it suffices, by the Grad Theorem, to prove
that \verb|eqweqmap| has a homotopy inverse.  To construct the
homotopy inverse we employ induction on weak equivalences.  That
is, to construct a map from the space \verb|weq A B| to the path
space \verb|paths A B| it suffices to be able to specify the
image of an identity weak equivalence.  But these are clearly sent to the
identity path.  The fact that this determines a homotopy inverse is
then a consequence of the computation principle for weak
equivalences.  (Full details of the proof can be found in the Coq file
accompanying this paper.)

\subsection{Function extensionality}

Henceforth, we assume that the Univalence Axiom holds.  That is, we
adopt the following:
\begin{center}
  \begin{coqcode}
Axiom univ : isweqeqweqmap.
  \end{coqcode}
\end{center}
In Coq, the command \verb|Axiom| serves to introduce a new global
hypothesis.  In this case, we assume given a proof \verb|univ| of
the Univalence Axiom.  Note that, although a good number do, not all
of the facts we prove below require the Univalence Axiom for their
proofs.  We will begin by briefly summarizing one of Voevodsky's type
theoretic results regarding the Univalence Axiom.

One somewhat curious feature of type theory without the Univalence
Axiom is that the principle of Function Extensionality is not
derivable.  Although it can be formulated in a number of ways, the
principle of Function Extensionality should be understood as stating
that, for continuous maps $f,g\colon A\to B$, paths from $f$ to $g$ in the
function space $B^{A}$ correspond to homotopies from $f$ to $g$.

Voevodsky \cite{Vo2012a} showed that Function Extensionality (and a
number of closely related principles) is a consequence of the
Univalence Axiom.  A sketch of the proof written in the usual
mathematical style can be found in Gambino \cite{Gambino:OR}.  We will
not describe the proof of Function Extensionality here.  We merely
mention that we will make use of it in what follows.  Explicitly, we
make use of a term 
\begin{center}
  \begin{coqcode}
funextsec : forall B : UU, forall E : B -> UU, 
 forall f g : ( forall x : B, E x ), 
  isweq ( pathtohtpysec f g ).
  \end{coqcode}
\end{center}
This is a slightly more general form of Function Extensionality and
implies the more common version as follows:
\begin{center}
  \begin{coqcode}
Definition funextfun { A B : UU } ( f g : A -> B ) : 
homot f g -> paths f g 
 := weqinv ( pair _ ( funextsec A ( fun z => B ) f g ) ).
  \end{coqcode}
\end{center}

\subsection{Impredicativity of h-levels}

We will now describe a number of consequences of the Univalence Axiom
concerning h-levels.  First, we give the definition of h-levels as follows:
\begin{center}
  \begin{coqcode}
Fixpoint isofhlevel ( n : nat ) ( A : UU ) := 
  match n with
    | O => iscontr A 
    | S n => forall a b : A, isofhlevel n ( paths a b )
  end.
  \end{coqcode}
\end{center}
Here the operation \verb|Fixpoint| tells Coq that we will define
a function out of an inductive type (in this case \verb|nat|)
recursively.  So, the definition is saying that a type $A$ is of
h-level 0 when it is contractible and it is of h-level $(n+1)$ when,
for all points $a$ and $b$ of $A$, the path space
$\texttt{paths}\;a\;b$ is of h-level $n$.

The h-levels satisfy an important property which type theorists refer
to as being \emph{impredicative}: they are closed under universal
quantification in the sense described below.  In the base case (for contractible
spaces), this is proved as follows:
\begin{center}
  \begin{coqcode}
Lemma impredbase { B : UU } ( E : B -> UU ) : 
(forall x : B, iscontr ( E x )) -> iscontr (forall x : B, E x).
  \end{coqcode}
\end{center}
Intuitively, what this says is that given a fibration $E$ over $B$, if
every fiber $E_{x}$ is contractible, then the space 
\verb|forall x : B, E x | of sections of the fibration is also
contractible.  We omit the proof, which is an immediate consequence of
\verb|funextsec|.

The following general principle of impredicativity of h-levels then
follows by induction:
\begin{center}
  \begin{coqcode}
Lemma impred ( n : nat ) : forall B : UU, forall E : B -> UU, 
( forall x : B, isofhlevel n ( E x ) ) 
    -> isofhlevel n ( forall x : B, E x ).
  \end{coqcode}
\end{center}
Again, this lemma states that if all fibers $E_{x}$ are of h-level
$n$, then so is the space of sections of the fibration.

\subsection{The total space and h-levels}

Next, we would like to explain the behavior of h-levels when it comes
to forming the total space of a fibration.  Assume given a fibration
$E\to B$ over $B$.  That is, we assume given a term 
\verb|E : B -> UU |.  Then, for any points $x$ and $y$, there is
a weak equivalence between the path space $\texttt{paths}\;x\;y$ and
the space which consists of pairs $(f,g)$ consisting of paths $f$ from $\pi_{1}(x)$ to
$\pi_{1}(y)$ and paths $g$ from $f_{!}(\pi_{2}(x))$ to $\pi_{2}(y)$
(see Section \ref{sec:transport} above for more on this idea).  Using this fact we obtain the
following lemma:
\begin{center}
  \begin{coqcode}
Lemma totalandhlevel ( n : nat ) : forall B : UU, 
forall E : B -> UU, forall is : isofhlevel n B, 
forall is' : ( forall x : B, isofhlevel n ( E x ) ), 
isofhlevel n ( total E ).   
  \end{coqcode}
\end{center}
The lemma states that if the base space $B$ and all fibers $E_{x}$ are
of h-level $n$, then so is the total space.  In the base case $n=0$
this is straightforward.  In the induction case, we observe that, for
any $x$ and $y$ in the total space,
\begin{center}
  \begin{coqcode}
isofhlevel n ( paths x y )
  \end{coqcode}
\end{center}
can be replaced, by the Univalence Axiom and the weak equivalence
mentioned above, by 
\begin{center}
  \begin{coqcode}
isofhlevel n ( 
total ( fun v : paths ( pr1 x ) ( pr1 y ) => 
paths ( transportf E v ( pr2 x ) ) ( pr2 y ) ) ).
  \end{coqcode}
\end{center}
But the induction hypothesis applies in this case, since we are now
dealing with a space of the form \verb|total ...|, and so we are done.

\subsection{The unit type and contractibility}

The unit type \verb|unit| corresponds to the terminal object $1$
in the category of spaces under consideration.  It is the inductive
type with a single generator \verb|tt : unit|.  For any type 
\verb|A| there is an induced map
\begin{center}
  \begin{coqcode}
tounit A : A -> unit.
  \end{coqcode}
\end{center}
It is a useful fact about contractible spaces that \verb|A| is
contractible if and only if \verb|tounit A| is a weak
equivalence.  We omit the straightforward proof of this equivalence.  

One fact about contractible spaces we will require is the fact that if
\verb|A| is contractible, then so is the type 
\verb|iscontr A| of proofs that \verb|A| contractible. This
is captured by the following lemma:
\begin{center}
  \begin{coqcode}
Lemma iscontrcontr { A : UU } ( is : iscontr A ) : 
iscontr ( iscontr A ).    
  \end{coqcode}
\end{center}
By the Univalence Axiom and the characterization of contractible
spaces mentioned above, in order to prove this theorem it suffices to
consider the case where \verb|A| is the unit type, which is
more or less immediate.  This proof reveals one important method for
using the Univalence Axiom: to prove something about a space $A$ it
suffices, by the Univalence Axiom, to prove the same fact about an
easier to manage space which is weakly equivalent to $A$.

\subsection{Some propositions}

We think of types of h-level 1 as being \emph{propositions} (or
truth-values) in the sense familiar to logicians.  Following this
intuition, we introduce the following notation:
\begin{center}
  \begin{coqcode}
Notation isaprop := ( isofhlevel 1 ).
  \end{coqcode}
\end{center}
Being a proposition is the same as being \emph{proof irrelevant}.
That is, $P$ is a proposition if and only if, for all terms $p,q:P$, there is
a path from $p$ to $q$.

One important consequence of \verb|iscontrcontr| is the fact that
being contractible is itself a proposition:
\begin{center}
  \begin{coqcode}
Lemma isapropiscontr ( A : UU ) : isaprop ( iscontr A ).
  \end{coqcode}
\end{center}
First, note that it is clear that a sufficient condition for being a
proposition is being contractible.  So, given points $p$ and $q$ of
type \verb|iscontr A|, it suffices to show that the type 
\verb|iscontr A| is itself contractible, which is by
\verb|iscontrcontr|.

As a consequence of impredicativity of h-levels together with
\verb|isapropiscontr| we obtain the following lemma:
\begin{center}
  \begin{coqcode}
Lemma isapropisweq { A B : UU } ( f : A -> B ) : isaprop ( isweq f ).
  \end{coqcode}
\end{center}
Again using impredicativity of h-levels and \verb|isapropisweq|
we obtain the following theorem:
\begin{center}
  \begin{coqcode}
Theorem isaxiomunivalence : isaprop ( isweqeqweqmap ).
  \end{coqcode}
\end{center}
That is, the type of the Univalence Axiom is a proposition and
therefore, assuming that there exists a term \verb|univ| of this
type, the space of such terms is contractible.

Similarly, a straightforward argument shows that, for any space $A$,
the type \verb|isofhlevel n A| is a proposition.

\subsection{The h-levels of h-universes}

We will now consider types of the form
\begin{center}
  \begin{coqcode}
total ( fun x : UU => isofhlevel n x )
  \end{coqcode}
\end{center}
which correspond to the types of all small spaces of a fixed
h-level $n$.  That is, they are what you might call
\emph{h-universes}.  Ultimately we will compute the h-levels of
h-universes.  First we will develop some further basic facts
about h-levels.

Note that if $A$ is of h-level $n$ then a straightforward argument
(using the discussion of h-levels of total spaces above), shows that
$A$ is also of h-level $(n+1)$.  That is, the h-universes are
cumulative.  Next, observe that if $B$ is of h-level $(n+1)$, then so
is the space of weak equivalences $\texttt{WEq}\;A\;B$ for any space
$A$:
\begin{center}
  \begin{coqcode}
Lemma hlevelweqcodomain ( n : nat ) : forall A B : UU, 
isofhlevel ( S n ) B -> isofhlevel ( S n ) ( weq A B ).
Proof.
  intros A B is. apply totalandhlevel. apply impred. 
  intros. apply is. intros. apply isaproptoisofhlevelSn. 
  apply isapropisweq.
Defined.
  \end{coqcode}
\end{center}
The proof is as follows.  It suffices, by 
\verb|totalandhlevel|, to prove separately that both the function
space \verb|A -> B| and the space of proofs that an element $f$ of
the function space is a weak equivalence are of h-level $(n+1)$.  The
former is by impredicativity of h-levels and the latter is by 
\verb|isapropisweq|.

Now, for each fixed $n$, there is a version of \verb|eqweqmap|
relativized to the h-universe of type of h-level $n$:
\begin{center}
  \begin{coqcode}
Definition hlevelneqweqmap ( n : nat ) 
( A B : total ( fun x : UU => isofhlevel n x ) ) 
: paths A B -> weq ( pr1 A ) ( pr1 B ) := 
fun f => eqweqmap ( pathintotalfiberpr1 f ) .
  \end{coqcode}
\end{center}
It turns out that, because \verb|isweq f| is a proposition, this
map is a weak equivalence:
\begin{center}
  \begin{coqcode}
Lemma isweqhlevelneqweqmap ( n : nat ) : 
forall A B : total ( fun x : UU => isofhlevel n x ), 
isweq ( hlevelneqweqmap n A B ).
  \end{coqcode}
\end{center}
Next, we observe that if both $A$ and $B$ are contractible, then so is
the space of weak equivalences from $A$ to $B$:
\begin{center}
\begin{coqcode}
Lemma iscontrandweq { A B : UU } ( is : iscontr A ) 
 ( is' : iscontr B ) : iscontr ( weq A B ).
\end{coqcode}
\end{center}
Finally, we observe that the h-universe of types of h-level $n$ itself
has h-level $(n+1)$:
\begin{center}
  \begin{coqcode}
Theorem isofhlevelSnhn ( n : nat ) : isofhlevel ( S n ) 
( total ( fun x : UU => isofhlevel n x ) ).
  \end{coqcode}
\end{center}
The proof is as follows.  First, note that if we are given $A$ and $B$
of type 
\begin{center}
  \begin{coqcode}
total ( fun x : UU => isofhlevel n x ),
  \end{coqcode}
\end{center}
then these terms should themselves be thought of as pairs $A=(A,p)$ and
$B=(B,q)$ where $p$ is a proof that $A$ is of h-level $n$ and $q$ is a
proof that $B$ is of h-level $n$.  Nonetheless, there is a weak
equivalence between the type $\texttt{paths}\;(A,p)\;(B,q)$ and the
type $\texttt{WEq}\;A\;B$ and as such it suffices to construct a
term of type
\begin{center}
  \begin{coqcode}
isofhlevel n ( weq A B ).
  \end{coqcode}
\end{center}
When $n=0$ this is by \verb|iscontrandweq| and in the case where
$n=m+1$ it is by \verb|hlevelweqcodomain|.

\emph{Acknowledgements.} We thank Pierre-Louis Curien, Marcelo Fiore, Helmut Hofer and Vladimir
Voevodsky for useful discussions on the topic of this paper.  We
also thank Chris Kapulkin for several helpful comments on a
preliminary draft.

\bibliographystyle{amsplain}
\bibliography{ref}

\end{document}